\documentclass[reqno,12pt]{amsart}
\usepackage[colorlinks=true, linkcolor=blue, citecolor=blue]{hyperref}

\usepackage{amssymb}
\usepackage{amsmath, graphicx, rotating}
\usepackage{color}
\usepackage{soul}
\usepackage[dvipsnames]{xcolor}

\allowdisplaybreaks
\usepackage{ifthen}
\usepackage{xkeyval}
\usepackage{todonotes}
\usepackage[T1]{fontenc}
\usepackage{lmodern}
\usepackage[english]{babel}

\usepackage{ upgreek }
\usepackage{stmaryrd}
\SetSymbolFont{stmry}{bold}{U}{stmry}{m}{n}
\usepackage{amsthm}
\usepackage{float}

\usepackage{ bbm }
\usepackage{ stmaryrd }
\usepackage{ mathrsfs }
\usepackage{ frcursive }
\usepackage{ comment }

\usepackage{pgf, tikz}
\usetikzlibrary{shapes}
\usepackage{varioref}
\usepackage{enumitem}

\definecolor{rouge}{rgb}{0.7,0.00,0.00}
\definecolor{vert}{rgb}{0.00,0.5,0.00}
\definecolor{bleu}{rgb}{0.00,0.00,0.8}
\usepackage[margin=1in]{geometry}
\newtheorem{theorem}{Theorem}[section]
\newtheorem*{theorem*}{Theorem}
\newtheorem{lemma}[theorem]{Lemma}

\newtheorem{proposition}[theorem]{Proposition}

\labelformat{hypothesis}{\textbf{M\kern-0.1mm#1}}

\newtheorem{condition}{Condition}

\newtheorem{conditionA}{A\kern-0.1mm}
\labelformat{conditionA}{\textbf{A\kern-0.1mm#1}}

\theoremstyle{definition}

\newtheorem{remark}[theorem]{Remark}

\def \eref#1{\hbox{(\ref{#1})}}

\numberwithin{equation}{section}

\def\geq{\geqslant}
\def\leq{\leqslant}

\def\RR{\mathbb{R}}
\def\PP{\mathbb{P}}
\def\EE{\mathbb{E}}

\def\vare{{\varepsilon}}
\def \eref#1{\hbox{(\ref{#1})}}

\def\EE{\mathbb{ E}}

\begin{document}

\title[Strong and weak convergence for multiscale SPDEs ]
{Orders of strong and weak averaging principle for multiscale SPDEs driven by $\alpha$-stable process}

\author{Xiaobin Sun}
\curraddr[Sun, X.]{ School of Mathematics and Statistics, Research Institute of Mathematical Science, Jiangsu Normal University, Xuzhou, 221116, China}
\email{xbsun@jsnu.edu.cn}

\author{Yingchao Xie}
\curraddr[Xie, Y.]{ School of Mathematics and Statistics, Research Institute of Mathematical Science, Jiangsu Normal University, Xuzhou, 221116, China}
\email{ycxie@jsnu.edu.cn}

\begin{abstract}
In this paper, the averaging principle is studied for a class of multiscale stochastic partial differential equations driven by $\alpha$-stable process, where $\alpha\in(1,2)$.
Using the technique of Poisson equation, the orders of strong and weak convergence are given $1-1/\alpha$ and $1-r$ for any $r\in (0,1)$ respectively.
The main results extend Wiener noise considered by Br\'{e}hier in \cite{B3} and Ge et al. in \cite{GSX} to $\alpha$-stable process, and the finite dimensional case considered by Sun et al.
in \cite{SXX} to the infinite dimensional case.
\end{abstract}

\date{\today}

\maketitle

\section{Introduction}

Many systems change involving slow and fast components in the natural world. For instance, dynamics of chemical reaction networks often take place on notably different times
scales, from the order of nanoseconds ($10^{-9}$ s) to the order of several days; When you are looking at the interaction between temperature and climate, it is found that
the daily temperature changes more rapidly, while climate changes are relatively slow. People always call this kind of system as the multiscale system or slow-fast system.
Multiscale models have wide applications in various fields, such as nonlinear oscillations, chemical kinetics, biology, climate dynamics, see e.g. \cite{BR, WTRY16} and the references therein.

\vspace{0.2cm}
Multiscale systems often show characteristics that do not conform to common sense, and the complexity of this kind system makes the traditional single theory no longer applicable,
so the study of multiscale models system becomes inevitable and necessary. The mathematical methods people use are often referred to as the methods of averaging and of homogenization,
see e.g. \cite{FW,PS} and the references therein.

\vspace{0.2cm}
The averaging principle for multiscale models describes the asymptotic behavior of the slow component as the scale parameter $\vare\to 0$. Bogoliubov and
Mitropolsky \cite{BM} first studied the averaging principle for the deterministic systems. Khasminskii \cite{K1} established an averaging principle for the stochastic differential
equations driven by Wiener noise. Since these pioneering works, many people have studied averaging principles for various stochastic systems, see e.g. \cite{GKK,GJ,HL,KY2,Ki,L1,LRSX1,RSX2,V0,ZHWWD}
for stochastic differential equations (SDEs), and see e.g. \cite{B1, B3, C1,C2,CL,DSXZ,FLL,FWLL,GP4,GP5,PXY,SXX1,WR,XML1} for stochastic partial differential equations (SPDEs).

\vspace{0.2cm}
In this paper, we consider the following slow-fast stochastic system on a Hilbert space $H$:
\begin{equation}\left\{\begin{array}{l}\label{main equation}
\displaystyle
dX^{\vare}_t=\left[AX^{\vare}_t+B(X^{\vare}_t, Y^{\vare}_t)\right]dt+dL_t,\quad X^{\vare}_0=x\in H,\vspace{2mm}\\
\displaystyle dY^{\vare}_t=\frac{1}{\vare}[AY^{\vare}_t+F(X^{\vare}_t, Y^{\vare}_t)]dt+\frac{1}{\vare^{1/\alpha}}dZ_t,\quad Y^{\vare}_0=y\in H,\end{array}\right.
\end{equation}
where $\varepsilon >0$ is a small parameter describing the ratio of time scales between the slow component $X^{\varepsilon}$ and fast component $Y^{\varepsilon}$.
$A$ is a selfadjoint operator, measurable functions $B,F:H\times H\rightarrow H$ satisfy some appropriate conditions, and $\{L_t\}_{t\geq 0}$ and $\{Z_t\}_{t\geq 0}$ are
mutually independent cylindrical $\alpha$-stable process with $\alpha\in (1,2)$, which are defined on a complete filtered probability space
$(\Omega,\mathscr{F},\{\mathscr{F}_{t}\}_{t\geq0},\mathbb{P})$.

\vspace{2mm}
The strong averaging principle for such stochastic systems \eref{main equation} has attracted some attention recently. For instance, Bao et al. \cite{BYY} proved
the strong averaging principle for two-time scale SPDEs driven by $\alpha$-stable noise. The authors have proved the strong averaging principle for stochastic Ginzburg-Landau equation,
stochastic Burgers equations and a class of SPDEs with H\"{o}lder coefficients in \cite{SZ}, \cite{CSS} and \cite{SXXZ} respectively. However, the key technique used in these mentioned papers
was based on the Khasminskii's time discretization, thus no satisfactory convergence order was obtained. Meanwhile, studying the convergence rate is an interesting and important topic
in multiscale system. For instance, Br\'{e}hier \cite{B2,B3} used the convergence rate to construct the efficient numerical schemes, based on the Heterogeneous Multiscale Methods.

\vspace{2mm}
The order of convergence for slow-fast stochastic systems has been studied extensively.
The technique of Khasminskii's time discretization is frequently used to study the strong convergence rate (see e.g. \cite{B1, GD, L1, RSX2}),
while the method of asymptotic expansion of solutions of Kolmogorov equations in the parameter $\vare$ is  used to study the  weak convergence rate (see e.g. \cite{B1,DSXZ, FWLL,KY2,ZFWL}).
Recently, the technique of Poisson equation is widely used to study the strong and weak convergence rates, see e.g. \cite{B3,GSX,RSX2,RXY,SXX,XY}. For more applications of Poisson equation,
see e. g. \cite{PV1,PV2,RX} and references therein.

\vspace{2mm}
The aim of this paper is first establish the strong convergence rates of stochastic system \eref{main equation}. More precisely, for any $(x,y)\in H^{\eta}\times H$ with $\eta\in (0,1)$, $T>0$
and $1\leq p<\alpha$, one tries to prove that
\begin{eqnarray*}
\sup_{t\in [0, T]}\mathbb{E}|X^{\vare}_t-\bar{X}_t|^p\leq C\vare^{\left(1-\frac{1}{\alpha}\right)p},
\end{eqnarray*}
where $C$ is a constant depending on $T, \|x\|_{\eta}, |y|, p$ and $\bar{X}$ is the solution of the corresponding averaged equation  (see Eq. \eref{1.3} below).

Secondly, we continuous to study the weak convergence rates of stochastic system \eref{main equation}. More precisely, for some fixed test function $\phi$, then for any $(x,y)\in H \times H$,
$T>0$ and $r\in (0,1)$, one tries to prove that
\begin{eqnarray*}
|\mathbb{E}\phi(X^{\vare}_t)-\EE\phi(\bar{X}_t)|\leq C\vare^{1-r},
\end{eqnarray*}
where $C$ is a constant depending on $T, |x|, |y|, r$.

\vspace{2mm}
In contrast to the existing works \cite{B3, GSX}, due to the Wiener noise is considered there, thus the solution has finite second moment usually.
However the solution here does not has finite second moment due to the $\alpha$-stable noise, hence some methods developed there do not work in this situation.
In order to overcome this difficulty, we shall estimate the solution of the corresponding Poisson equation more carefully, meanwhile the accurate treatment of the
$\alpha$-stable process is provided in the proof.

\vspace{2mm}
In contrast to the existing work \cite{SXX}, the strong and weak convergence rates for slow-fast SDEs driven by $\alpha$-stable noise are obtained there,
we here extend the case of finite dimension to infinite dimension essentially. However, we have to overcome some non-trivial difficulties in the infinite dimensional case.
For example the presentation of term $AX^{\vare}_t$, the method of Galerkin approximation and the smoothing properties of the semigroup $e^{tA}$ will be used to deal with
a serious of difficulties arising from the unbounded operator $A$.

\vspace{2mm}
Another contribution of this paper is to fill a gap in \cite{BYY} partially. As stated in \cite[Remark 3.3]{BYY}, \emph{"for the technical reason, it seems hard to show
Theorem 3.1 without the uniform boundedness of the nonlinearity"}, where \emph{"Theorem 3.1"} means the strong averaging principle holds and \emph{"the nonlinearity"} means the coefficient $B$.
In fact, the essential reason is that the method used in \cite{BYY} is the classical Khasminskii's time discretization, which highly depends on the square calculation in the
proof, hence the finite second moment of the solution $X^{\vare}_t$ is required usually. But the solution $X^{\vare}_t$ for system \eref{main equation} only has finite $p$-th moment
($0<p<\alpha$), the uniform boundedness of $B$ is used to weaken the finite second moment to finite first moment. However, the technique of Poisson equation is used to remove
the condition of uniform boundedness of $B$, but some bounded conditions of second and third  derivatives for the coefficients are assumed. Moreover,
the optimal strong averaging convergence order is obtained here.

\vspace{2mm}
The organization of this paper is as follows. In the next section, some notations and assumptions are introduced. Then we state our main results.
Section 3 is devoted to study the regularity of the solution of the corresponding Poisson equation. The detailed proofs of strong and weak convergence rates are provided
in Sections 4 and 5 respectively. The final section is the appendix, where we give some a-priori estimates of the solution, and study the Galerkin approximation of the system
\eref{main equation} and the finite dimensional approximation of the frozen equation.

\vspace{0.2cm}
We note that throughout this paper $C$, $C_p$, $C_{T}$ and $C_{p,T}$ denote positive constants which may change from line to line, where the subscript $p,T$
are used to emphasize that the constants depend on $p,T$.

\section{Notations and main results} \label{Sec Main Result}

\subsection{Notations and assumptions}

We introduce some notation used throughout this paper.  $H$ is Hilbert space with inner product $\langle\cdot,\cdot\rangle$ and norm $|\cdot|$.
$\mathbb{N}_{+}$ stands for the collection of all the positive integers.

\vspace{2mm}
$\mathcal{B}(H)$ denotes the collection of all measurable functions $\varphi(x): H\rightarrow \RR$.
For any $k\in \mathbb{N}_{+}$,
\begin{align*}
C^k(H):=&\{\varphi\in \mathcal{B}(H): \varphi \mbox{ and all its Fr\'{e}chet derivatives up to order } k \mbox{ are continuous}\},\\
C^k_b(H):=&\{\varphi\in C^k(H): \mbox{ for } 1\le i\le k, \mbox{ all $i$-th Fr\'{e}chet derivatives of } \varphi \mbox{ are bounded} \}.
\end{align*}
For any $\varphi\in C^3(H)$, by the Riesz representation theorem, we often identify the first Fr\'{e}chet derivative $D\varphi(x)\in \mathcal{L}(H, \RR)\cong H$,
the second derivative $D^2\varphi(x)$ as a linear operator in $\mathcal{L}(H,H)$ and the third derivative $D^3\varphi(x)$ as a linear operator in $\mathcal{L}(H,\mathcal{L}(H,H))$, i.e.,
\begin{eqnarray*}
&&D\varphi(x) \cdot h_1 = \langle D\varphi(x), h_1 \rangle,\quad h_1\in H,\\
&&D^2\varphi(x) \cdot (h_1,h_2) = \langle D^2\varphi(x)\cdot h_1, h_2 \rangle,  \quad  h_1,h_2 \in H,\\
&&D^3\varphi(x) \cdot (h_1,h_2,h_3)= [D^3\varphi(x)\cdot h_1]\cdot (h_2,h_3),\quad h_1,h_2,h_3\in H,
\end{eqnarray*}
where $D^k\varphi(x) \cdot h$ is the $k$-th directional derivative of $\varphi$ in the direction $(h_1,\ldots,h_k)$, for $k=1,2,3$.
%
%

A selfadjoint operator $A$ satisfies $Ae_n=-\lambda_n e_n$ with $\lambda_n>0$ and $\lambda_n\uparrow \infty$, as $n\uparrow \infty$, where $\{e_n\}_{n\geq1}\subset \mathscr{D}(A)$
is a complete orthonormal basis of $H$. For any $s\in\RR$, we define
 $$H^s:=\mathscr{D}((-A)^{s/2}):=\left\{u=\sum_{k\in\mathbb{N}_{+}}u_ke_k: u_k\in \mathbb{R},~\sum_{k\in\mathbb{N}_{+}}\lambda_k^{s}u_k^2<\infty\right\}$$
and
 $$(-A)^{s/2}u:=\sum_{k\in \mathbb{N}_{+}}\lambda_k^{s/2} u_ke_k,~~u\in\mathscr{D}((-A)^{s/2})$$
with the associated norm $\|u\|_{s}:=|(-A)^{s/2}u|=\left(\sum_{k\in\mathbb{N}_{+}}\lambda_k^{s} u^2_k\right)^{1/2}$. It is easy to see $\|\cdot\|_0=|\cdot|$.

\vspace{2mm}
The following smoothing properties of the semigroup $e^{tA}$ (see \cite[Proposition 2.4]{B1}) will be used quite often later in this paper:
\begin{eqnarray}
\|e^{tA}x\|_{\sigma_2}\leq C_{\sigma_1,\sigma_2}t^{-\frac{\sigma_2-\sigma_1}{2}}e^{-\frac{\lambda_1 t}{2}}\|x\|_{\sigma_1},\quad x\in H^{\sigma_2},\sigma_1\leq\sigma_2, t>0,\label{P3}
\end{eqnarray}
\begin{eqnarray}
|e^{tA}x-x|\leq C_{\sigma}t^{\frac{\sigma}{2}}\|x\|_{\sigma},\quad x\in H^{\sigma},\sigma>0, t\geq 0.\label{P4}
\end{eqnarray}

Let $\{L_t\}_{t\geq 0}$ and $\{Z_t\}_{t\geq 0}$ be mutually independent cylindrical $\alpha$-stable processes, where $\alpha\in(1,2)$, i.e.,
$$
L_t=\sum_{k\in \mathbb{N}_{+}}\beta_{k}L^{k}_{t}e_k,\quad Z_t=\sum_{k\in \mathbb{N}_{+}}\gamma_{k}Z^{k}_{t}e_k,\quad t\geq 0,
$$
where $\{\beta_k\}_{k\in \mathbb{N}_{+}}$ and $\{\gamma_k\}_{k\in\mathbb{N}_{+}}$ are two given sequence of positive numbers, $\{L^n_t\}_{n\geq1}$ and $\{Z^n_t\}_{n\geq1}$
are two sequences of independent one dimensional rotationally symmetric $\alpha$-stable processes with $\alpha\in(1,2)$ satisfying for any $k\in \mathbb{N}_{+}$ and $t\geq0$,
$$\mathbb{E}[e^{i L^k_{t}h}]=\mathbb{E}[e^{i Z^k_{t}h}]=e^{-t|h|^{\alpha}}, \quad h\in \mathbb{R}.$$

For $t>0$, $k\in \mathbb{N}_{+}$ and $\Gamma\in\mathscr{B}(\mathbb{R}\setminus\{0\})$, the Poisson random measure associated with $L^{k}$ and $Z^k$ are defined by
$$N^{1,k}([0,t],\Gamma)=\sum_{0\leq s\leq t}1_{\Gamma}(L^{k}_s-L^{k}_{s-}),\quad N^{2,k}([0,t],\Gamma)=\sum_{0\leq s\leq t}1_{\Gamma}(Z^{k}_s-Z^{k}_{s-})$$
and the corresponding compensated Poisson random measures are given by
$$\widetilde{N}^{i,k}([0,t],\Gamma)=N^{i,k}([0,t],\Gamma)-t\nu(\Gamma),\quad i=1,2,$$
where $\nu(dy)=\frac{c_{\alpha}}{|y|^{1+\alpha}}dy$ is the L\'evy measure with $c_{\alpha}>0$.

\vspace{2mm}
By L\'evy-It\^o's decomposition and the symmetry of the L\'{e}vy measure $\nu$, one has
$$L^{k}_{t}=\int_{|x|\leq c}x\widetilde{N}^{1,k}([0,t],dx)+\int_{|x|>c}x N^{1,k}([0,t],dx),$$
$$Z^{k}_{t}=\int_{|x|\leq c}x\widetilde{N}^{2,k}([0,t],dx)+\int_{|x|>c}x N^{2,k}([0,t],dx),$$
where $c>0$. We also assume that $\{L^n_t\}_{n\geq1}$ and $\{Z^n_t\}_{n\geq1}$ are independent.

\vspace{0.2cm}

Now, we assume the following conditions on the coefficients $B, F: H\times H \rightarrow H$ throughout the paper:


\begin{conditionA}\label{A1}
$B$ and $F$ are Lipschitz continuous, i.e., there exist positive constants $L_{F}$ and $C$ such that for any $x_1,x_2,y_1,y_2\in H$,
\begin{eqnarray*}
&&\left|B(x_1,y_1)-B(x_2,y_2)\right|\leq C(|x_1-x_2|+|y_1-y_2|),\\
&&\left|F(x_1,y_2)-F(x_2,y_2)\right|\leq C|x_1-x_2|+L_{F}|y_1-y_2|.
\end{eqnarray*}
\end{conditionA}

\begin{conditionA}\label{A2} Assume that $\lambda_{1}-L_{F}>0$, $\sum_{k\in \mathbb{N}_{+}}\beta^{\alpha}_k\lambda^{\alpha-1}_k<\infty$ and $\sum_{k\in \mathbb{N}_{+}}\gamma^{\alpha}_k<\infty$.
\end{conditionA}

\begin{conditionA}\label{A3}
Assume that there exists $\kappa_1\in (0,2)$ such that the following directional derivatives are well-defined and satisfy:
\begin{equation}\left\{\begin{array}{l}\label{aa}
\displaystyle
|D_{x}B(x,y)\cdot h|\leq C |h| \quad \text{and} \quad |D_{y}B(x,y)\cdot h|\leq C |h|,\quad \forall x,y, h\in H ,\vspace{0.2cm} \\
|D_{xx}B(x,y)\cdot(h,k)|\leq C |h|\|k\|_{\kappa_1},\quad \forall x,y, h\in H, k\in H^{\kappa_1},\vspace{0.2cm} \\
|D_{yy}B(x,y)\cdot(h,k)|\leq C |h|\|k\|_{\kappa_1},\quad \forall x,y, h\in H, k\in H^{\kappa_1} ,\vspace{0.2cm}\\
|D_{xy}B(x,y)\cdot (h,k)|\leq C |h|\|k\|_{\kappa_1}, \quad \forall x,y, h\in H, k\in H^{\kappa_1},\vspace{0.2cm}\\
|D_{xy}B(x,y)\cdot (h,k)|\leq C \|h\|_{\kappa_1}|k|, \quad \forall x,y, k\in H, h\in H^{\kappa_1},\vspace{0.2cm}\\
|D_{xyy}\!B(x,y)\cdot(h,k,l)|\leq C|h|\|k\|_{\kappa_1}\|l\|_{\kappa_1}, \quad \forall x,y, h\in H, k,l\in H^{\kappa_1},\vspace{0.2cm}\\
|D_{yyy}\!B(x,y)\cdot(h,k,l)|\leq C|h|\|k\|_{\kappa_1}\|l\|_{\kappa_1}, \quad \forall x,y, h\in H, k,l\in H^{\kappa_1} ,\vspace{0.2cm}\\
|D_{xxy}\!B(x,y)\cdot(h,k,l)|\leq C|h|\|k\|_{\kappa_1}\|l\|_{\kappa_1}, \quad \forall x,y, h\in H, k,l\in H^{\kappa_1} ,
\end{array}\right.
\end{equation}
where $D_{x}B(x,y)\cdot h$ is the directional derivative of $B(x,y)$ in the direction $h$ with respective to $x$,
other notations can be interpreted similarly. The above properties in \eref{aa} also hold for operator $F$.
\end{conditionA}

\begin{remark}\label{Re1} Suppose that the assumption \ref{A1} holds, $\sum_{k\in \mathbb{N}_{+}}\beta^{\alpha}_k/\lambda_k<\infty$ and
$\sum_{k\in \mathbb{N}_{+}}\gamma^{\alpha}_k/{\lambda_k}<\infty$. As \cite{PZ} did, we have that, for $\varepsilon>0$ and $(x,y)\in H\times H$,
the system \eref{main equation} admits a unique mild solution $(X^{\vare}_t, Y^{\vare}_t)\in H\times H$, i.e., $\PP$-a.s.,
\begin{equation}\left\{\begin{array}{l}\label{A mild solution}
\displaystyle
X^{\varepsilon}_t=e^{tA}x+\int^t_0e^{(t-s)A}B(X^{\varepsilon}_s, Y^{\varepsilon}_s)ds+\int^t_0 e^{(t-s)A}dL_s,\vspace{0.2cm}\\
\displaystyle
Y^{\varepsilon}_t=e^{tA/\varepsilon}y+\frac{1}{\varepsilon}\int^t_0e^{(t-s)A/\varepsilon}F(X^{\varepsilon}_s,Y^{\varepsilon}_s)ds
+\frac{1}{\vare^{1/\alpha}}\int^t_0 e^{(t-s)A/\varepsilon}dZ_s.
\end{array}\right.
\end{equation}
\end{remark}

\begin{remark}\label{Re2}
The condition $\lambda_{1}-L_{F}>0$ in assumption \ref{A2} is called the strong dissipative condition, which is used to prove the existence and uniqueness of the invariant measures
and the exponential ergodicity of the transition semigroup of the frozen equation. The condition
$\sum_{k\in \mathbb{N}_{+}}\gamma^{\alpha}_k<\infty$ in assumption \ref{A2} are necessary when applying It\^o's formula for the solution $(X^{\varepsilon}, Y^{\varepsilon})$.
While the condition $\sum_{k\in \mathbb{N}_{+}}\beta^{\alpha}_k \lambda^{\alpha-1}_k<\infty$ is used to control $\EE\|\int^t_0 e^{(t-s)A}d L_s\|_2$.
For a more general result see \cite[Lemma 4.1]{PSXZ}, i.e., if $\sum_{k\in \mathbb{N}_{+}}\frac{\beta^{\alpha}_k}{\lambda^{1-\alpha\theta/2}_k}<\infty$ holds for some
$\theta\geq 0$, then we have for any  $0<p<\alpha$,
\begin{eqnarray}
\sup_{t\geq 0}\EE\left\|\int^t_0 e^{(t-s)A}d L_s\right\|^p_{\theta}\leq C_{\alpha,p}\left(\sum_{k\in\mathbb{N}_{+}}\frac{\beta^{\alpha}_k}{\lambda^{1-\alpha\theta/2}_k}\right)^{p/\alpha}.\label{LA}
\end{eqnarray}
\end{remark}

\begin{remark} \label{EA3}
Here we give a example that the conditions in assumption \ref{A3} hold.
Let $H:=\{g\in L^2(\mathcal{D}): g(\xi)=0, \xi\in \partial\mathcal{D}\}$, where $\mathcal{D}=(0,1)^d$ for $d=1,2,3$, $\partial\mathcal{D}$ is the boundary of domain $\mathcal{D}$ and
$A:=\Delta:=\sum^d_{i=1}\partial^2_{\xi_i}$ be the Laplacian operator. The coefficient $B$ is defined to be the Nemytskii operator associated with a function
$b:\RR\times \RR\rightarrow \RR$, i.e., $B(x,y)(\xi):=b(x(\xi),y(\xi))$. Then the following directional derivatives are well-defined and belong to $H$,
\begin{eqnarray*}
&&D_{x}B(x,y)\cdot h=\partial_x b(x,y)h \quad \text{and} \quad D_{x}B(x,y)\cdot h=\partial_y b(x,y)h ,\quad \forall x,y, h\in H ,\vspace{2mm}\\
&&D_{xx}B(x,y)\cdot(h,k)=\partial_{xx}b(x,y)hk,\quad \forall x,y, h\in H, k\in L^{\infty}(0,1) ,\vspace{2mm}\\
&&D_{yy}B(x,y)\cdot(h,k)=\partial_{yy}b(x,y)hk,\quad \forall x,y, h\in H, k\in L^{\infty}(0,1),\vspace{2mm}\\
&&D_{xy}B(x,y)\cdot (h,k)=\partial_{xy}b(x,y)hk, \quad \forall x,y, h\in H, k\in L^{\infty}(0,1) ,\vspace{2mm}\\
&&D_{xyy}B(x,y)\cdot(h,k,l)=\partial_{xyy}b(x,y)hkl, \quad \forall x,y, h\in H, k,l\in L^{\infty}(0,1) ,\vspace{2mm}\\
&&D_{yyy}B(x,y)\cdot(h,k,l)=\partial_{yyy}b(x,y)hkl, \quad \forall x,y, h\in H, k,l\in L^{\infty}(0,1) ,\vspace{2mm}\\
&&D_{xxy}B(x,y)\cdot(h,k,l)=\partial_{xxy}b(x,y)hkl, \quad \forall x,y, h\in H, k,l\in L^{\infty}(0,1),
\end{eqnarray*}
where all the partial derivatives of $b$ appear above are uniformly bounded by assumption.
\vspace{2mm}

Note that $H^{d/2+\eta}\subset L^{\infty} (\mathcal{D})$ for any $\eta>0$(see \cite[(6)]{B3}),  thus it is easy to see that $B$ satisfy assumption \ref{A3}
with any $\kappa_1\in (1/2,2)$ for $d=1$, $\kappa_1\in (1,2)$ for $d=2$ and $\kappa_1\in (3/2,2)$ for $d=3$.
\end{remark}

\subsection{Main results}

Let $\mu^{x}$ be the unique invariant measure of the transition semigroup of the frozen equation
\begin{eqnarray}\label{FZE}
\left\{ \begin{aligned}
&dY_{t}=\left[AY_{t}+F(x,Y_{t})\right]dt+d Z_{t},\vspace{2mm}\\
&Y_{0}=y\in H
\end{aligned} \right.
\end{eqnarray}
and define $\bar{B}(x):=\int_{H}B(x,y)\mu^{x}(dy)$. Let $\bar{X}$ be the solution of the corresponding averaged equation:
\begin{equation}\left\{\begin{array}{l}
\displaystyle d\bar{X}_{t}=\left[A\bar{X}_{t}+\bar{B}(\bar{X}_{t})\right]dt+d L_{t},\vspace{2mm}\\
\bar{X}_{0}=x\in H. \end{array}\right. \label{1.3}
\end{equation}
Here we state our main results.
\begin{theorem} (\textbf{Strong convergence rate})\label{main result 1}
Suppose that assumptions \ref{A1} and \ref{A3} hold. Then for any initial values $(x,y)\in H^{\eta}\times H$  with $\eta\in(0,1)$, $T>0$, $1\leq p<\alpha$ and small enough $\vare,\delta>0$, we have
\begin{align}
\sup_{t\in [0,T]} \mathbb{E}|X_{t}^{\vare}-\bar{X}_{t}|^{p}\leq C_{p,T,\delta}\left[1+\|x\|^{(1+\delta)p}_{\eta}+|y|^{(1+\delta)p}\right]\vare^{\left(1-\frac{1}{\alpha}\right)p}. \label{ST}
\end{align}
\end{theorem}
\begin{remark}
The result \eref{ST} above implies that the strong convergence order is $1-\frac{1}{\alpha}$, which is the optimal order in the strong sense usually (see \cite[Example 2.2]{SXX}).
 Meanwhile, when $\alpha\uparrow2$, this order $1-\frac{1}{\alpha}\uparrow \frac{1}{2}$, which is in accord with the optimal order $1/2$ in the case of Wiener noise (see \cite{B3,GSX,RXY}).
 Note that we do not assume the boundedness of $B$, thus it gives a positive answer to \cite[Remark 3.3]{BYY}.
\end{remark}

\begin{theorem}(\textbf{Weak convergence rate})\label{main result 2}
Suppose that assumptions \ref{A1}-\ref{A3} hold. Moreover, $\sup_{x,y}|B(x,y)|<\infty$ and there exists $\kappa_2\in (0,2)$ such that the following directional derivatives are well-defined and satisfy:
\begin{eqnarray}
&&\!\!\!\!\!\!\!\!\!\!\!\!\!\!\!\!\!|D_{xxx}B(x,y)\cdot(h,k,l)|\leq C|h|\|k\|_{\kappa_2}\|l\|_{\kappa_2}, \quad \forall x,y, h\in H, k,l\in H^{\kappa_2}, \label{ThirdDer F4}\vspace{2mm}\\
&&\!\!\!\!\!\!\!\!\!\!\!\!\!\!\!\!\!|D_{xxx}F(x,y)\cdot(h,k,l)|\leq C|h|\|k\|_{\kappa_2}\|l\|_{\kappa_2}, \quad \forall x,y, h\in H, k,l\in H^{\kappa_2}. \label{ThirdDer F5}
\end{eqnarray}
Then for any test function $\phi\in C^3_b(H)$, initial values $(x,y)\in H\times H$, $T>0$, $r\in (0,1)$ and $\vare>0$, we have
\begin{align}
\sup_{t\in [0,T]}\left|\mathbb{E}\phi(X_{t}^{\vare})-\EE\phi(\bar{X}_{t})\right|\leq C_{r,T,\delta}\left[1+|x|^{1+\delta}+|y|^{1+\delta}\right]\vare^{1-r}, \label{WC}
\end{align}
where $C_{r,T,\delta}$ is a constant depends on $r, T, \delta$ and $\lim_{r\downarrow 0}C_{r,T,\delta}=\infty$.
\end{theorem}

\begin{remark}
The result \eref{WC} implies that the weak convergence rate is $1-r$ with any $r\in (0,1)$. Comparing with the theorem \ref{main result 1}, the stronger regularity of the coefficients
$B,F$ are assumed, while initial value $x\in H$ and the improved convergence order is obtained. It is worthy to point that the boundedness of the $B$ is assumed for the reason of the
solution does not has finite second moment. Meanwhile, it fails to obtain the expected weak convergence order 1 (see \cite{SXX}).
\end{remark}
\begin{remark}
Here we give a example that some additional conditions in Theorem \eref{main result 2} hold. Recall the notations in Remark \ref{EA3}. Obviously $\sup_{x,y\in H}|B(x,y)|<\infty$ by
assuming $\sup_{x,y\in \RR}|b(x,y)|<\infty$. Next we check the condition \eref{ThirdDer F5}. Assume that $\sup_{x,y\in \RR}|\partial_{xxx}b(x,y)|<\infty$. Then the following directional
derivatives are well-defined and belong to $H$,
$$
D_{xxx}B(x,y)\cdot(h,k,l)=\partial_{xxx}b(x,y)hkl,\quad \forall x,y, h\in H, k,l\in L^{\infty}(0,1).
$$
Note that $H^{d/2+\eta}\subset L^{\infty} (\mathcal{D})$ for any $\eta>0$,  thus it is easy to see that
$$
|D_{xxx}B(x,y)\cdot(h,k,l)|=|\partial_{xxx}b(x,y)hkl|\leq C|h|\|k\|_{\kappa_2}\|l\|_{\kappa_2},
$$
where $\kappa_2\in (1/2,2)$ for $d=1$, $\kappa_2\in (1,2)$ for $d=2$ and $\kappa_2\in (3/2,2)$ for $d=3$. The assumption \eref{ThirdDer F5} can be handled similarly.
\end{remark}

\section{The Poisson equation for nonlocal operator}

Since the drift coefficient $B$ may not be bounded and the solution $X^{\vare}_t$ does not has finite second moment, the classical Khasminskii's time discreatization dose not work in
this situation (see \cite[Remark 3.3]{BYY}). We shall use the technique of Poisson equation to obtain the strong and weak convergence rates for system \eref{main equation}.
Meanwhile, note that the operator $A$ is not a bounded operator and $H\ni x$ may not belong to $\mathscr{D}(-A)$, we use Galerkin approximation to reduce the infinite dimensional
problem to a finite dimension firstly, then we will take the limit finally, i.e., considering
\begin{equation}\left\{\begin{array}{l}\label{Ga mainE}
\displaystyle
dX^{m,\vare}_t=[AX^{m,\vare}_t+B^m(X^{m,\vare}_t, Y^{m,\vare}_t)]dt+d\bar L^m_t,\  X^{m,\vare}_0=x^{m}\in H_m,\vspace{2mm}\\
\displaystyle
dY^{m,\vare}_t=\frac{1}{\vare}[AY^{m,\vare}_t+F^m(X^{m,\vare}_t, Y^{m,\vare}_t)]dt+\frac{1}{\vare^{1/\alpha}}d\bar Z^m_t,\quad Y^{m,\vare}_0=y^m\in H_m,\end{array}\right.
\end{equation}
where $m\in \mathbb{N}_{+}$, $H_{m}:=\text{span}\{e_{k};1\leq k \leq m\}$, $\pi_{m}$ is the orthogonal projection of $H$ onto $H_{m}$, $x^{m}:=\pi_m x, y^m:=\pi_m y$ and
\begin{eqnarray*}
&& B^{m}(x,y):=\pi_m B(x,y),\quad F^{m}(x,y):=\pi_m F(x, y),\\
&&\bar L^{m}_t:=\sum^m_{k=1}\beta_{k}L^{k}_{t}e_k,\quad \bar Z^{m}_t:=\sum^m_{k=1}\gamma_{k}Z^{k}_{t}e_k.
\end{eqnarray*}

Similarly, we consider the following approximation to the averaged equation \eref{1.3}:
\begin{equation}\left\{\begin{array}{l}
\displaystyle d\bar{X}^m_{t}=\left[A\bar{X}^m_{t}+\bar{B}^m(\bar{X}^m_{t})\right]dt+d \bar{L}^m_{t},\vspace{2mm}\\
\bar{X}^m_{0}=x^m,\end{array}\right. \label{Ga 1.3}
\end{equation}
where $\bar{B}^m(x):=\int_{H_m}B^m(x,y)\mu^{x,m}(dy)$, and $\mu^{x,m}$ is the unique invariant measure of the transition semigroup of the following frozen equation:
\begin{equation}\left\{\begin{array}{l}
\displaystyle dY^{x,y,m}_t=[AY^{x,y,m}_t+F^m(x,Y^{x,y,m}_t)]dt+d\bar{Z}^m_t,\vspace{2mm}\\
Y^{x,y,m}_0=y\in H_m.\end{array}\right. \label{Ga FZE}
\end{equation}

Before we use the technique of Poisson equation, we need to do some preparations.


\subsection{The frozen equation}
For fixed $x\in H$ and $m\in\mathbb{N}_{+}$, we recall the finite dimensional frozen equation \eref{Ga FZE}.
Note that $F^m(x,\cdot)$ is Lipschitz continuous, then it is easy to show that for any initial value $y\in H_m$,
equation $(\ref{Ga FZE})$ has a unique mild solution $\{Y_{t}^{x,y,m}\}_{t\geq 0}$ in $H_m$, i.e., $\PP$-a.s.,
\begin{eqnarray*}
Y_{t}^{x,y,m}=e^{tA}y+\int^t_0 e^{(t-s)A}F^m(x,Y_{s}^{x,y,m})ds+\int^t_0e^{(t-s)A}d\bar{Z}^m_s.
\end{eqnarray*}
Moreover, the solution $\{Y_{t}^{x,y,m}\}_{t\geq 0}$ is a time homogeneous Markov process. Let $\{P^{x,m}_t\}_{t\geq 0}$ be its transition semigroup, i.e., for any bounded
measurable function $\varphi: H_m\rightarrow \mathbb{R}$,
$$
P^{x,m}_t\varphi(y):=\EE\varphi(Y_{t}^{x,y,m}), \quad y\in H_m, t\geq 0.
$$

Before studying the asymptotic behavior of $\{P^{x,m}_t\}_{t\geq 0}$, we prove the following lemmas.

\begin{lemma}\label{L3.2}
For any $p\in [1,\alpha)$, there exists $C_p>0$ such that

\begin{eqnarray}
\sup_{m\geq 1}\EE|Y_{t}^{x,y,m}|^p\leq e^{-\lambda_1 pt}|y|^p+C_p(1+|x|^p).\label{FEq2}
\end{eqnarray}
\end{lemma}
\begin{proof}
By a straightforward computation, it is easy to see
\begin{eqnarray*}
|Y^{x,y,m}_t|\leq\!\!\!\!\!\!\!\!&&e^{-\lambda_1 t}|y|+\int^t_0e^{-\lambda_1(t-s)}|F^m(x,Y^{x,y,m}_s)|ds+\left|\int^t_0 e^{(t-s)A}dZ_s\right|\\
\leq\!\!\!\!\!\!\!\!&& e^{-\lambda_1 t}|y|+\int^t_0e^{-\lambda_1(t-s)}\left(L_F|Y^{x,y,m}_s|+C|x|+C\right)ds+\left|\int^t_0 e^{(t-s)A}dZ_s\right|.
\end{eqnarray*}
Then by Minkowski's inequality, we have for any $p\in (1,\alpha)$ and $t\geq 0$,
\begin{eqnarray*}
\left(\EE|Y^{x,y,m}_t|^p\right)^{1/p}\leq\!\!\!\!\!\!\!\!&&e^{-\lambda_1 t}|y|+\int^t_0e^{-\lambda_1(t-s)}\left[L_F\left(\EE|Y^{x,y,m}_s|^p\right)^{1/p}+C|x|+C\right]ds\nonumber\\
&&+\left(\EE\left|\int^t_0 e^{(t-s)A}dZ_s\right|^p\right)^{1/p}\\
\leq\!\!\!\!\!\!\!\!&& e^{-\lambda_1 t}|y|+\frac{C(1+|x|)}{\lambda_1}+\frac{L_F}{\lambda_1}\sup_{s\in[0,t]}\left(\EE|Y^{x,y,m}_s|^p\right)^{1/p}+\left(\EE\left|\int^t_0 e^{(t-s)A}dZ_s\right|^p\right)^{1/p}.
\end{eqnarray*}
By condition $L_F<\lambda_1$ in assumption \ref{A2}, we obtain
\begin{eqnarray}
\EE|Y^{x,y,m}_t|^p
\leq\!\!\!\!\!\!\!\!&& e^{-\lambda_1 pt}|y|^p+C_p(1+|x|^p)+C_p\sup_{t\geq 0}\EE\left|\int^t_0 e^{(t-s)A}dZ_s\right|^p.\label{F3.13}
\end{eqnarray}

Note that condition $\sum^{\infty}_{k=1}\gamma^{\alpha}_k<\infty$ implies $\sum^{\infty}_{k=1}\frac{\gamma^{\alpha}_k }{\lambda_k}<\infty$. Then by \eref{LA}, we get
\begin{eqnarray}
\sup_{t\geq 0}\EE\left|\int^t_0 e^{(t-s)A}dZ_s\right|^p
\leq C_{\alpha,p}\left(\sum^{\infty}_{k=1}\frac{\gamma^{\alpha}_k }{\lambda_k} \right)^{p/\alpha}.\label{F3.14}
\end{eqnarray}

Hence, \eref{F3.13} and \eref{F3.14} implies that \eref{FEq2} holds. The proof is complete.
\end{proof}

\begin{lemma}\label{L3.6}
For any $t>0$ and $x_1,x_2\in H, y_1,y_2\in H_m$, there exists $C>0$ such that
\begin{eqnarray}
\sup_{m\geq 1}|Y^{x_1,y_1,m}_t-Y^{x_2,y_2,m}_t|\leq e^{-\frac{(\lambda_1-L_F) t}{2}}|y_1-y_2|+C|x_1-x_2|.\label{increase Y}
\end{eqnarray}
\end{lemma}
\begin{proof}
For any $x_1,x_2\in H, y_1,y_2\in H_m$, note that
\begin{eqnarray*}
\frac{d}{dt}(Y^{x_1,y_1,m}_t-Y^{x_2,y_2,m}_t)=\!\!\!\!\!\!\!\!&&A(Y^{x_1,y_1,m}_t-Y^{x_2,y_2,m}_t)+\left[F^m(x_1, Y^{x_1,,y_1,m}_t)-F^m(x_2, Y^{x_2,y_2,m}_t)\right].
\end{eqnarray*}
Multiplying $2(Y^{x_1,y_1,m}_t-Y^{x_2,y_2,m}_t)$ in both sides, we obtain
\begin{eqnarray*}
\frac{d}{dt}|Y^{x_1,y_1,m}_t-Y^{x_2,y_2,m}_t|^2=\!\!\!\!\!\!\!\!&&-2\|Y^{x_1,y_1,m}_t-Y^{x_2,y_2,m}_t\|^2_1\\
&&+2\langle F^m( x_1,Y^{x_1,y_1,m}_t)-F^m(x_2, Y^{x_2,y_2,m}_t), Y^{x_1,y_1,m}_t-Y^{x_2,y_2,m}_t\rangle.
\end{eqnarray*}
Then by $L_F<\lambda_1$ in assumption \ref{A2} and Young's inequality, we get
\begin{eqnarray*}
\frac{d}{dt}|Y^{x_1,y_1,m}_t-Y^{x_2,y_2,m}_t|^2\leq\!\!\!\!\!\!\!\!&&-2\lambda_1|Y^{x_1,y_1,m}_t-Y^{x_2,y_2,m}_t|^2+2L_F\left|Y^{x_1,y_1,m}_t-Y^{x_2,y_2,m}_t\right|^2\\
&&+C|x_1-x_2|\left|Y^{x_1,y_1,m}_t-Y^{x_2,y_2,m}_t\right|\nonumber\\
\leq\!\!\!\!\!\!\!\!&& -(\lambda_1-L_F)\left|Y^{x_1,y_1,m}_t-Y^{x_2,y_2,m}_t\right|^2+C|x_1-x_2|^2.
\end{eqnarray*}
By comparison theorem, we have for any $t> 0$,
\begin{eqnarray}
|Y^{x_1,y_1,m}_t-Y^{x_2,y_2,m}_t|^2\leq e^{-(\lambda_1-L_F)t}|y_1-y_2|^2+C|x_1-x_2|^2.\label{increase Y^2}
\end{eqnarray}
The proof is complete.
\end{proof}

\vspace{0.2cm}
Under the condition $\lambda_1-L_F>0$, it is well known that the transition semigroup $\{P^{x,m}_t\}_{t\geq 0}$ admits a unique invariant measure $\mu^{x,m}$ (see e.g.
\cite[Theorem 1.1]{WJ}). Using \eref{FEq2}, it is easy to check that for any $p\in [1,\alpha)$ and $x\in H$,
\begin{eqnarray}
\sup_{m\geq 1}\int_{H_m}|z|^p\mu^{x,m}(dz)\leq C_p(1+|x|^p).\label{E3.20}
\end{eqnarray}
Furthermore, we shall prove the following exponential ergodicity for the transition semigroup $\{P^{x,m}_t\}_{t\geq 0}$,  which plays an important role in studying the regularity of the
solution of the Poisson equation.
\begin{proposition} \label{ergodicity in finite}
For any Lipschitz continuous function $G: H\rightarrow H$, then we have for any $t>0$,
\begin{eqnarray}
\sup_{x\in H,m\geq 1}\left| P^{x,m}_tG(y)-\mu^{x,m}(G)\right|\leq\!\!\!\!\!\!\!\!&& C\|G\|_{Lip} e^{-\frac{(\lambda_1-L_F) t}{2}}(1+|x|+|y|), \label{ergodicity1}
\end{eqnarray}
where $C>0$ and $\|G\|_{Lip}:=\sup_{x\neq y\in H}\frac{|G(x)-G(y)|}{|x-y|}$.
%
\end{proposition}
\begin{proof}
By the definition of invariant measure $\mu^{x,m}$ and \eref{increase Y},  we have for any $t> 0$,
\begin{eqnarray*}
\left| P^{x,m}_tG(y)-\mu^{x,m}(G)\right|=\!\!\!\!\!\!\!\!&&\left| \EE G( Y^{x,y,m}_t)-\int_{H_m}G(z)\mu^{x,m}(dz)\right|\\
\leq\!\!\!\!\!\!\!\!&& \left|\int_{H_m}\left[\EE G(Y^{x,y,m}_t)-\EE G(Y^{x,z,m}_t)\right]\mu^{x,m}(dz)\right|\\
\leq\!\!\!\!\!\!\!\!&& \|G\|_{Lip}\int_{H_m} \EE\left| Y^{x,y,m}_t-Y^{x,z,m}_t\right|\mu^{x,m}(dz)\\
\leq\!\!\!\!\!\!\!\!&& \|G\|_{Lip} e^{-\frac{(\lambda_1-L_F) t}{2}}\int_{H_m}|y-z|\mu^{x,m}(dz)\\
\leq\!\!\!\!\!\!\!\!&& \|G\|_{Lip}e^{-\frac{(\lambda_1-L_F) t}{2}}\left[|y|+\int_{H_m}|z|\mu^{x,m}(dz)\right]\\
\leq\!\!\!\!\!\!\!\!&& C \|G\|_{Lip}e^{-\frac{(\lambda_1-L_F) t}{2}}(1+|x|+|y|),
\end{eqnarray*}
where the last inequality is a consequence of \eref{E3.20}. The proof is complete.


\end{proof}


\subsection{The regularity of solution of the Poisson equation}
By the preparation in above subsection,  this subsection is devoted to study the following Poisson equation:
\begin{equation}
-\mathscr{L}^m_{2}(x)\Phi_m(x,y)=B^m(x,y)-\bar{B}^m(x),\quad x,y\in H_m,\label{PE2}
\end{equation}
where $\mathscr{L}^m_{2}(x)$ is the infinitesimal generator of the transition semigroup of the finite dimensional frozen equation \eref{Ga FZE}, i.e.,
\begin{eqnarray}
&&\mathscr{L}^m_{2}(x)\Phi_m(x,y)\nonumber\\
=\!\!\!\!\!\!\!\!&&D_y\Phi_m(x,y)\cdot (Ay+F^m(x,y))\nonumber\\
\!\!\!\!\!\!\!\!&&+\sum^m_{k=1}\gamma^{\alpha}_k \int_{\RR}\!\!\Phi_m(x,y+e_k z)\!-\!\Phi_m(x,y)\!-\!\langle D_y\Phi_m(x,y), e_k z\rangle 1_{\{|z|\leq 1\}}\nu (dz).\label{L_2}
\end{eqnarray}
Note that $\mathscr{L}_{2}(x)$ is the infinitesimal generator of the transition semigroup of the frozen process $\{Y^{x,y,m}_t\}_{t\geq 0}$, we define
\begin{eqnarray}
\Phi_m(x,y):=\int^{\infty}_{0}\left[\EE B^m(x,Y^{x,y,m}_t)-\bar{B}^m(x)\right]dt.\label{SPE}
\end{eqnarray}
It is easy to check \eref{SPE} solves equation \eref{PE2}. The following is the regularity of the solution $\Phi_m(x,y)$ with respect to parameters, which will play an important
role in the proof of our main results. The regularity of the solution of the Poisson equation with respect to parameters have been study in some references, see e.g. \cite{PV2,RSX2,SXX}.


\begin{proposition}\label{P3.6}
For any $\delta\in (0,1]$, there exists $C, C_{\delta}>0$ such that for any $x,y, h,k\in H_m$,
\begin{eqnarray}
&&\sup_{m\geq1}|\Phi_m(x,y)|\leq C(1+|x|+|y|); \label{E1}\\
&&\sup_{m\geq 1}|D_y \Phi_m(x,y)\cdot h|\leq C|h|;\label{E21}\\
&&\sup_{m\geq 1}|D_x \Phi_m(x,y)\cdot h|\leq C_{\delta}(1+|x|^{\delta}+|y|^{\delta})|h|;\label{E2}\\
&& \sup_{m\geq 1}| D_{xx}\Phi_m(x, y)\cdot(h,k)|\leq C_{\delta}(1+|x|^{\delta}+|y|^{\delta})|h|\|k\|_{\kappa_1},\label{E3}
\end{eqnarray}
where $\kappa_1$ is the constant in assumption \ref{A3}.
\end{proposition}
\begin{remark}
It is crucial that the parameter $\delta>0$ in \eref{E2} and \eref{E3}, which plays a key role in solving the difficulty cased by the solution does not has finite second moment.
\end{remark}
\begin{proof}
The proof is divided into three steps.

\vspace{0.2cm}
\textbf{Step 1.}
By Proposition \ref{ergodicity in finite}, we get
\begin{eqnarray*}
|\Phi_m(x,y)|\leq\!\!\!\!\!\!\!\!&&\int^{\infty}_{0}| \EE B^m(x,Y^{x,y,m}_t)-\bar{B}^m(x)|dt\\
\leq\!\!\!\!\!\!\!\!&& C(1+|y|)\int^{\infty}_{0}e^{-\frac{(\lambda_1-L_F) t}{2}}dt\leq C(1+|y|).
\end{eqnarray*}
So the first estimate in \eref{E1} holds.

For any $h\in H_m$, we have
$$D_y \Phi_m(x,y)\cdot h=\int^{\infty}_0 \EE[D_y B^m(x,Y^{x,y,m}_t)\cdot D^h_y Y^{x,y,m}_t]dt,$$
where $D^h_{y} Y^{x,y,m}_t$ is the derivative of $Y^{x,y,m}_t$ with respect to $y$ in the direction $h$, which satisfies
 \begin{equation}\left\{\begin{array}{l}\label{partial y}
\displaystyle
dD^h_{y} Y^{x,y,m}_t=AD^h_{y} Y^{x,y,m}_tdt+D_y F^m(x,Y^{x,y}_{t})\cdot D^h_{y} Y^{x,y,m}_tdt,\vspace{2mm}\\
D^h_{y} Y^{x,y}_0=h.\\
\end{array}\right.
\end{equation}
Then by $\lambda_1-L_F>0$ and condition \eref{aa},  it is easy to see
\begin{eqnarray}
\sup_{x,y\in H_m}|D^h_y Y^{x,y,m}_t|\leq Ce^{-\frac{(\lambda_1-L_F)t}{2}}|h|,\label{partial yY}
\end{eqnarray}
and for any $\eta\in (0,2)$ and $t>0$,
\begin{eqnarray}
\sup_{x,y\in H_m}\|D^h_y Y^{x,y,m}_t\|_{\eta}\leq Ce^{-\frac{(\lambda_1-L_F)t}{2}}(1+t^{-\eta/2})|h|.\label{Hpartial yY}
\end{eqnarray}
Thus it follows
$$\sup_{x,y\in H_m}|D_y \Phi_m(x,y)\cdot h|\leq C|h|.$$
So  \eref{E2} holds.

Now, we define
\begin{eqnarray*}
\tilde B^m_{t_0}(x, y, t):=\hat{B}^m(x,y, t)-\hat{B}^m(x, y, t+t_0),
\end{eqnarray*}
where $\hat{B}^m(x, y, t):=\EE B^m(x, Y^{x, y}_t)$. Note that \eref{ergodicity1} implies
$$
\lim_{t_0\rightarrow \infty} \tilde{B}^m_{t_0}(x, y, t)=\EE[B^m(x,Y^{x,y,m}_t)]-\bar{B}^m(x).
$$
So in order to prove \eref{E2} and \eref{E3},  it suffices to show that there exists $C>0$ such that for any $\delta\in (0,1]$, $t_0>0$, $ t\geq 0$, $x,y\in H$,
\begin{eqnarray}
|D_{x}\tilde{B}^m_{t_0}(x,y,t)\cdot h|\leq Ce^{-\frac{(\lambda_1-L_F)\delta t}{4}}(1+t^{-\frac{\kappa_1 \delta}{2}})(1+|x|^{\delta}+|y|^{\delta})|h|,\label{E21}
\end{eqnarray}
\begin{eqnarray}
\left|D_{xx} \tilde{B}^m_{t_0}(x,y,t)\cdot(h,k)\right|\leq Ce^{-\frac{(\lambda_1-L_F)\delta t}{4}}(1+t^{-\frac{\kappa_1 \delta}{2}})(1+|x|^{\delta}+|y|^{\delta})|h|\|k\|_{\kappa_1}\label{E22},
\end{eqnarray}
which will be proved in step 2 and step 3 respectively.
\vspace{0.3cm}

\textbf{Step 2.} In this step, we intend to prove \eref{E21}. It follows from the Markov property,
\begin{eqnarray*}
\tilde{B}^m_{t_0}(x,y, t)=\!\!\!\!\!\!\!\!&& \hat{B}^m(x, y, t)-\EE B^m(x, Y^{x,y,m}_{t+t_0})\nonumber\\
=\!\!\!\!\!\!\!\!&& \hat{B}(x,y, t)-\EE \{\EE[B^m(x,Y^{x,y,m}_{t+t_0})|\mathscr{F}_{t_0}]\}\nonumber\\
=\!\!\!\!\!\!\!\!&& \hat{B}(x, y, t)-\EE \hat{B}^m(x, Y^{x,y,m}_{t_0},t).
\end{eqnarray*}
Then we obtain
\begin{eqnarray}
D_{x}\tilde{B}^m_{t_0}(x,y, t)\cdot h=\!\!\!\!\!\!\!\!&& D_{x} \hat{B}^m(x, y, t)\cdot h-\EE D_{x}\hat{B}^m(x, Y^{x,y,m}_{t_0},t)\cdot h\nonumber\\
&&- \EE \left[D_y\hat{B}^m(x,Y^{x,y,m}_{t_0},t)\cdot D^h_{x} Y^{x,y,m}_{t_0} \right],\label{5.8}
\end{eqnarray}
where $D^h_{x} Y^{x,y,m}_t$ is the derivative of $Y^{x,y,m}_t$ with respect to $x$ in the direction $h$, which satisfies
 \begin{equation*}\left\{\begin{array}{l}\label{V}
\displaystyle
dD^h_{x} Y^{x,y,m}_t=\left[AD^h_{x} Y^{x,y,m}_t+D_{x} F^m(x,Y^{x,y,m}_{t})\cdot h+D_y F^m(x,Y^{x,y,m}_{t})\cdot D^h_{x} Y^{x,y,m}_t\right]dt,\vspace{2mm}\\
D^h_{x} Y^{x,y,m}_0=0.
\end{array}\right.
\end{equation*}
By  \eref{P3} and \eref{aa}, we can easily obtain for any $\eta\in [0,2)$
\begin{eqnarray}
\sup_{t\geq 0, x,y\in H_m} \|D^h_{x} Y^{x,y,m}_t\|_{\eta}\leq C|h|.\label{HS0}
\end{eqnarray}

Note that
$$D_y\hat{B}^m(x, y, t)\cdot h=\EE\left[D_y B^m(x, Y^{x,y,m}_t)\cdot D^h_y Y^{x,y,m}_t\right],$$
$$D_x\hat{B}^m(x, y, t)\cdot h=\EE\left[D_x B^m(x, Y^{x,y,m}_t)\cdot h+ D_y B^m(x, Y^{x,y,m}_t)\cdot D^h_x Y^{x,y,m}_t\right].$$
Then by \eref{aa}, \eref{HS0} and \eref{partial yY}, we get
\begin{eqnarray}
\sup_{x, y\in H_m}|D_y\hat{B}^m(x, y, t)\cdot h|\leq Ce^{-\frac{(\lambda_1-L_F)t}{2}}|h|,\label{S1}
\end{eqnarray}
\begin{eqnarray}
\sup_{t\geq 0, x, y\in H_m}|D_x\hat{B}^m(x, y, t)\cdot h|\leq  C|h|.\label{S11}
\end{eqnarray}

Next if we can prove that for any $t> 0, h,k\in H_m$,
\begin{equation}\label{rs}
\sup_{x,y\in H_m}|D_{xy}\hat B^m(x,y, t)\cdot (h,k)|\leq C e^{-\frac{(\lambda_1-L_F) t}{4}}(1+t^{-\frac{\kappa_1}{2}})|h||k|,
\end{equation}
which combins with \eref{S11} and \eref{FEq2} we obtain for any $\delta\in (0,1]$,
\begin{eqnarray}
&&|D_{x} \hat B^m(x, y, t)\cdot h-\EE D_{x}\hat B^m(x, Y^{x,y,m}_{t_0},t)\cdot h|\nonumber\\
\leq\!\!\!\!\!\!\!\!&&C|h|^{1-\delta}|D_{x} \hat B^m(x, y, t)\cdot h-\EE D_{x}\hat B^m(x, Y^{x,y,m}_{t_0},t)\cdot h|^{\delta}\nonumber\\
\leq\!\!\!\!\!\!\!\!&& C|h|^{1-\delta}\left[\EE\left |\int^1_0 D_{xy}\hat B^m( x, \xi y+(1-\xi)Y^{x,y,m}_{t_0}, t)\cdot (h, y-Y^{x,y,m}_{t_0})d\xi\right|\right]^{\delta}\nonumber\\
\leq\!\!\!\!\!\!\!\!&& Ce^{-\frac{(\lambda_1-L_F)\delta t}{4}}(1+t^{-\frac{\kappa_1 \delta}{2}})\EE|y-Y^{x,y,m}_{t_0}|^{\delta}|h|\nonumber\\
\leq\!\!\!\!\!\!\!\!&& Ce^{-\frac{(\lambda_1-L_F)\delta t}{4}}(1+t^{-\frac{\kappa_1 \delta}{2}})(1+|x|^{\delta}+|y|^{\delta})|h|.\label{S2}
\end{eqnarray}

Thus by \eref{S1} and \eref{S2}, we get
\begin{eqnarray*}
|D_{x}\tilde B^m_{t_0}( x, y, t)\cdot h|\leq\!\!\!\!\!\!\!\!&& Ce^{-\frac{(\lambda_1-L_F)\delta t}{2}}(1+|x|^{\delta}+|y|^{\delta})|h|+Ce^{-\frac{(\lambda_1-L_F)t}{2}}|h|\\
\leq\!\!\!\!\!\!\!\!&&Ce^{-\frac{(\lambda_1-L_F)\delta t}{2}}(1+|x|^{\delta}+|y|^{\delta})|h|,
\end{eqnarray*}
which proves \eref{E21}.


Now, we are in a position to prove \eref{rs}. Note that
\begin{eqnarray*}
&&|D_{xy}\hat B^m(x,y, t)\cdot (h,k)|=|D_{xy} \left[\EE B^m(x, Y^{x,y,m}_t)\right]\cdot (h,k)|\nonumber\\
=\!\!\!\!\!\!\!\!&& |D_y\left[\EE D_{x} B^m(x, Y^{x,y,m}_t)+\EE D_{y} B^m(x, Y^{x,y,m}_t)\cdot D_{x} Y^{x,y,m}_t\right]\cdot (h,k)|\\
\leq\!\!\!\!\!\!\!\!&& \left|\EE \left[D_{xy} B^m(x, Y^{x,y,m}_t)\cdot(D^h_{y}Y^{x,y,m}_t,k)\right]\right|+\left|\EE \left[D_{y} B^m(x, Y^{x,y,m}_t)\cdot D^{(h,k)}_{xy} Y^{x,y,m}_t\right]\right|\\
&&+\left|\EE \left[D_{yy} B^m(x, Y^{x,y,m}_t)\cdot(D^h_{x} Y^{x,y,m}_t, D^k_{y} Y^{x,y,m}_t)\right]\right|,
\end{eqnarray*}
where $D^{(h,k)}_{xy}Y^{x,y,m}_t$ is the second
order derivative of $Y^{x,y,m}_{t}$ in the direction $(h,k)\in H_m\times H_m$ (once with respect to $x$ in the direction $h\in H_m$ and once with respect to $y$ in the direction $k\in H_m$),
which satisfies
 \begin{equation}\left\{\begin{array}{l}\label{Sxy DY}
\displaystyle
d D^{(h,k)}_{xy}Y^{x,y,m}_{t}=\left[AD^{(h,k)}_{xy}Y^{x,y,m}_{t}+D_{xy}F^m(x,Y^{x,y,m}_{t})\cdot (h,D^k_{y} Y^{x,y,m}_{t})\right.\nonumber\\
\displaystyle\quad\quad\quad\quad\quad\quad\quad+D_{yy}F^m(x,Y^{x,y,m}_{t})\cdot (D^{h}_{x}Y^{x,y,m}_{t},D^{k}_{y}Y^{x,y,m}_{t})\nonumber\\
\displaystyle\quad\quad\quad\quad\quad\quad\quad\left.+D_{y}F^m(x,Y^{x,y,m}_{t})\cdot D^{(h,k)}_{xy}Y^{x,y,m}_{t}\right]dt.\nonumber\\
\displaystyle D^{(h,k)}_{xy}Y^{x,y,m}_{0}=0.\\
\end{array}\right.
\end{equation}
By \eref{aa}, \eref{Hpartial yY},  \eref{HS0} and $\lambda_1-L_F>0$ , it is easy to prove for any $\eta\in [0,2)$ and $t>0$,
\begin{eqnarray}
\sup_{x,y\in H_m}\|D^{(h,k)}_{xy}Y^{x,y,m}_t\|_{\eta}\leq Ce^{-\frac{(\lambda_1-L_F)t}{4}}(1+t^{-\frac{\kappa_1}{2}})|h||k|.\label{Hpartial xyY}
\end{eqnarray}
Hence, \eref{rs} holds by assumption \ref{A3}, \eref{HS0}, \eref{Hpartial yY} and \eref{Hpartial xyY}.

\vspace{2mm}

\textbf{Step 3.} In this step, we intend to prove \eref{E22}. Recall that
\begin{eqnarray*}
D_{x}\tilde B^m_{t_0}(x,y, t)\cdot h=\!\!\!\!\!\!\!\!&& D_{x} \hat B^m(x, y, t)\cdot h-\EE D_{x}\hat B^m(x, Y^{x,y,m}_{t_0},t)\cdot h\\
&&- \EE \left[D_y\hat B^m(x,Y^{x,y,m}_{t_0},t)\cdot D^h_{x} Y^{x,y,m}_{t_0} \right].\label{5.8}
\end{eqnarray*}
Then it is easy to see
\begin{eqnarray*}
&&D_{xx}\tilde B^m_{t_0}(x,y, t)\cdot (h,k)\\
=\!\!\!\!\!\!\!\!&& \left[D_{xx} \hat B^m(x, y, t)\cdot(h,k)-\EE D_{xx}\hat B^m(x, Y^{x,y,m}_{t_0},t)\cdot(h,k)\right]\nonumber\\
&&-\EE \left[D_{xy}\hat B^m(x,Y^{x,y,m}_{t_0},t)\cdot (h, D^k_{x} Y^{x,y,m}_{t_0})+D_{yx}\hat B^m(x, Y^{x,y,m}_{t_0},t)\cdot(D^h_{x} Y^{x,y,m}_{t_0},k)\right]\\
&&-\EE\left[D_{yy}\hat B^m(x, Y^{x,y,m}_{t_0},t) \cdot (D^h_{x} Y^{x,y,m}_{t_0},D^k_{x} Y^{x,y,m}_{t_0})\right]\\
&&-\EE \left[ D_{y}\hat B^m(x,Y^{x,y,m}_{t_0},t) \cdot D^{(h,k)}_{xx} Y^{x,y,m}_{t_0} \right]:=\sum^4_{i=1} J_i,
\end{eqnarray*}
where $D^{(h,k)}_{xx}Y^{x,y,m}_{t}$ is the second derivative of $Y^{x,y,m}_{t}$
with respect to $x$ towards the direction $(h, k) \in H\times H$, which satisfies
 \begin{equation}\left\{\begin{array}{l}\label{SDY}
\displaystyle
d D^{(h,k)}_{xx}Y^{x,y,m}_{t}\!\!=\!\!\left[AD^{(h,k)}_{xx}Y^{x,y,m}_{t}\!+\!D_{xx}F^m(x,Y^{x,y,m}_{t})\cdot(h,k)\!+\!D_{xy}F^m(x,Y^{x,y,m}_{t})\cdot (h,D^k_{x} Y^{x,y,m}_{t})\right.\nonumber\\
\quad\quad\quad\quad+D_{yx}F^m(x,Y^{x,y,m}_{t})\cdot(D^h_{x} Y^{x,y,m}_{t},k)+D_{yy}F^m(x,Y^{x,y,m}_{t})\cdot (D^{h}_{x}Y^{x,y,m}_{t},D^{k}_{x}Y^{x,y,m}_{t})\nonumber\\
\quad\quad\quad\quad+\left.D_{y}F^m(x,Y^{x,y,m}_{t})\cdot D^{(h,k)}_{xx}Y^{x,y,m}_{t}\right]dt,\\
D^{(h,k)}_{xx}Y^{x,y,m}_{0}=0.\\
\end{array}\right.
\end{equation}

For the term $J_1$, note that
\begin{eqnarray*}
D_{x}\hat B^m(x,y, t)\cdot h=\!\!\!\!\!\!\!\!&&\EE \left[D_x B^m(x, Y^{x,y,m}_{t})\cdot h\right]+\EE\left[D_{y}B^m(x,Y^{x,y,m}_{t})\cdot D^h_{x}Y^{x,y,m}_t \right],
\end{eqnarray*}
which implies
\begin{eqnarray*}
&&D_{xx}\hat B^m(x,y, t)\cdot(h,k)\nonumber\\
=\!\!\!\!\!\!\!\!&& \EE \left[D_{xx}B^m(x,Y^{x,y,m}_t)\cdot(h,k)\right]+\EE \left[D_{xy} B^m(x,Y^{x,y,m}_t) \cdot (h,D^k_x Y^{x,y,m}_t)\right]\nonumber\\
&&+ \EE \left[D_{yx} B^m(x, Y^{x, y}_t)\cdot (D^h_{x} Y^{x,y,m}_{t},k)\right]+ \EE \left[D_{yy} B^m(x, Y^{x, y}_t)\cdot(D^h_{x} Y^{x,y,m}_{t}, D^k_{x} Y^{x,y,m}_{t})\right]\\
&&+ \EE \left[D_{y} B^m(x, Y^{x, y}_t)\cdot D^{(h,k)}_{xx}Y^{x,y,m}_{t} \right].
\end{eqnarray*}
Then it follows
\begin{eqnarray*}
&&D_{xxy}\hat B^m(x,y, t)\cdot(h,k,l)\\
=\!\!\!\!\!\!\!\!&& \EE \Big[D_{xxy}B^m(x,Y^{x,y,m}_t)\cdot(h,k,D^l_{y} Y^{x,y,m}_{t})+D_{xyy} B^m(x,Y^{x,y,m}_t) \cdot (h,D^k_x Y^{x,y,m}_t,D^l_{y} Y^{x,y,m}_{t})\nonumber\\
&&+D_{xy} B^m(x,Y^{x,y,m}_t) \cdot (h,D^{(k,l)}_{xy} Y^{x,y,m}_t)\!+\!D_{yxy} B^m(x, Y^{x,y,m}_t)\cdot (D^h_{x} Y^{x,y,m}_{t},k, D^l_{y} Y^{x,y,m}_{t})\\
&&+D_{yx} B^m(x,Y^{x,y,m}_t) \cdot (D^{(h,l)}_{xy} Y^{x,y,m}_t, k)\!+\!D_{yyy} B^m(x, Y^{x,y,m}_t)\cdot(D^h_{x} Y^{x,y,m}_{t}, D^k_{x} Y^{x,y,m}_{t},D^l_{y} Y^{x,y,m}_{t})\\
&&+D_{yy} B^m(x,Y^{x,y,m}_t) \cdot (D^{(h,l)}_{xy} Y^{x,y,m}_t, D^k_{x} Y^{x,y,m}_{t})\!+\!D_{yy} B^m(x,Y^{x,y,m}_t) \cdot (D^{h}_{x} Y^{x,y,m}_t, D^{(k,l)}_{xy} Y^{x,y,m}_{t})\\
&&+D_{yy} B^m(x, Y^{x,y,m}_t)\cdot (D^{(h,k)}_{xx}Y^{x,y,m}_{t}, D^l_{y} Y^{x,y,m}_{t})\!+\!D_{y} B^m(x, Y^{x,y,m}_t)\cdot D^{(h,k,l)}_{xxy}Y^{x,y,m}_{t} \Big],
\end{eqnarray*}
where $D^{(h,k,l)}_{xxy}Y^{x,y,m}_{t}$ is the
third order derivative of $Y^{x,y,m}_{t}$ (twice with respect to $x$ in the direction $(h,k)\in H_m\times H_m$ and once with respect to $y$ in the direction $l\in H_m$).

\vspace{2mm}
By assumption \ref{A3}, \eref{HS0} and \eref{Hpartial yY}, it is easy to prove that
\begin{eqnarray}
&&\sup_{t\geq 0,x,y\in H_m}|D^{(h,k)}_{xx} Y^{x,y,m}_t|\leq C|h||k|,\label{EDxx}
\end{eqnarray}
$$\sup_{ x,y\in H_m}|D^{(h,k,l)}_{xxy}Y^{x,y,m}_{t}|\leq Ce^{-\frac{(\lambda_1-L_F)t}{4}}(1+t^{-\frac{\kappa_1}{2}})|h||k||l|.$$
which combine with \eref{Hpartial yY} , \eref{HS0},  \eref{Hpartial xyY} and assumption \ref{A3} , we get
\begin{eqnarray}
|D_{xx}\hat B^m(x,y, t)\cdot(h,k)|\leq C|h|\|k\|_{\kappa_1},\label{DxxhatB}
\end{eqnarray}
\begin{eqnarray}
|D_{xxy}\hat B^m(x,y, t)\cdot(h,k,l)|\leq Ce^{-\frac{(\lambda_1-L_F)t}{4}}(1+t^{-\frac{\kappa_1}{2}})|h|\|k\|_{\kappa_1}|l|.\label{DxxyhatB}
\end{eqnarray}
Then by \eref{DxxhatB}, \eref{DxxyhatB} and similar as we did in proving \eref{S2}, we obtain for any $\delta\in (0,1]$,
\begin{eqnarray}
J_1\leq\!\!\!\!\!\!\!\!&& Ce^{-\frac{(\lambda_1-L_F)t}{4}}(1+t^{-\frac{\kappa_1}{2}})|h|\|k\|_{\kappa_1}\EE|y-Y^{x,y,m}_{t_0}|^{\delta}\nonumber\\
\leq\!\!\!\!\!\!\!\!&& Ce^{-\frac{(\lambda_1-L_F)t}{4}}(1+t^{-\frac{\kappa_1}{2}})|h|\|k\|_{\kappa_1}(1+|x|^{\delta}+|y|^{\delta}).\label{J1}
\end{eqnarray}

For the term $J_2$, note that
\begin{eqnarray*}
D_{xy}\hat B^m(x,y,t)\cdot(h,k)
=\!\!\!\!\!\!\!\!&& \EE \left[ D_{xy}B^m(x, Y^{x,y,m}_t)\cdot(h,D^k_{y}Y^{x,y,m}_t)\right]\\
&&+\EE \left[ D_{y}B^m(x,Y^{x,y,m}_t)\cdot D^{(h,k)}_{xy}Y^{x,y,m}_t\right]\\
&&+\EE\left[D_{yy}B^m(x,Y^{x,y,m}_{t})\cdot(D^h_{x}Y^{x,y,m}_t, D^k_{y}Y^{x,y,m}_t)\right].
\end{eqnarray*}
Then by assumption \ref{A3}, \eref{Hpartial yY}, \eref{HS0} and \eref{Hpartial xyY}, we have
$$
\sup_{x,y\in H}|D_{xy}\hat B^m(x,y,t)\cdot(h,k)|\leq C e^{-\frac{(\lambda_1-L_F)t}{4}}(1+t^{-\frac{\kappa_1}{2}})|h||k|.
$$
Thus it follows
\begin{eqnarray}
J_2\leq\!\!\!\!\!\!\!\!&& C e^{-\frac{(\lambda_1-L_F)t}{4}}(1+t^{-\frac{\kappa_1}{2}})|h||k|.\label{J2}
\end{eqnarray}

For the term $J_3$, by a similar argument as in estimating $J_2$, we have
$$
\sup_{x,y\in H_m}|D_{yy}\hat B^m(x,y,t)\cdot(h,k)|\leq C e^{-\frac{(\lambda_1-L_F)t}{4}}(1+t^{-\frac{\kappa_1}{2}})|h||k|.
$$
Hence, it is easy to see that
\begin{eqnarray}
J_3\leq C e^{-\frac{(\lambda_1-L_F)t}{4}}(1+t^{-\frac{\kappa_1}{2}})|h||k|.\label{J3}
\end{eqnarray}

For the term $J_4$, by \eref{EDxx} and \eref{S1},
we easily get
\begin{eqnarray}
J_4\leq C e^{-\frac{(\lambda_1-L_F)t}{2}}|h||k|.\label{J4}
\end{eqnarray}

Finally, \eref{E22} can be easily obtained by combining \eref{J1}-\eref{J4}. The proof is complete.
\end{proof}

\section{Strong convergence order}

This section is devoted to prove Theorem \ref{main result 1}. We first study the well-posedness of equation \eref{Ga 1.3} which approximates to the averaged equation.
Then we give the detailed proof of Theorem \ref{main result 1}.
\begin{lemma} \label{barX}
Equation \eref{Ga 1.3} exists a unique mild solution $\bar{X}^m_{t}$ satisfies
\begin{eqnarray}
\bar X^m_t=e^{tA}x^m+\int^t_0e^{(t-s)A} \bar B^m(\bar X^m_s)ds+\int^t_0 e^{(t-s)A}d\bar{L}^m_s.\label{3.22}
\end{eqnarray}
Moreover, for any $x\in H$, $T>0$ and $1\leq p< \alpha$, there exists a constant $C_{p,T}>0$ such that
\begin{align} \label{Control X}
\sup_{m\geq 1,t\in [0,T]}\mathbb{E}|\bar{X}^m_{t} |^p\leq  C_{p,T}(1+ |x|^p).
\end{align}
\end{lemma}
\begin{proof}
It is sufficient to check that the $\bar{B}^m$ is Lipschitz continuous, then \eref{Ga 1.3} admits a unique mild solution $\bar{X}^m_{t}$. The estimate \eref{Control X} can be proved by
the same argument as in the proof of Lemma \ref{GA1} in the appendix.

\vspace{2mm}
Indeed, for any $x_1,x_2\in H_m$ and $t>0$, by Proposition \ref{ergodicity1} and \eref{increase Y}, we have
\begin{eqnarray*}
|\bar{B}^m(x_1)-\bar{B}^m(x_2)|
=\!\!\!\!\!\!\!\!&&\left|\int_{H_m} B^m(x_1,y)\mu^{x_1,m}(dy)-\int_{H_m} B^m(x_2,y)\mu^{x_2,m}(dy)\right|\\
\leq\!\!\!\!\!\!\!\!&&\left|\EE B^m(x_1, Y^{x_1,0,m}_t)-\int_{H_m} B^m(x_1,z)\mu^{x_1,m}(dz)\right|\\
&&+\left|\EE B^m(x_2, Y^{x_2,0,m}_t)-\int_{H_m} B^m(x_2,z)\mu^{x_2,m}(dz)\right|\\
&&+\left|\EE B^m(x_1, Y^{x_1,0,m}_t)-\EE B^m(x_2, Y^{x_2,0,m}_t)\right|\\
\leq\!\!\!\!\!\!\!\!&&Ce^{-\frac{(\lambda_1-L_F) t}{2}}+C\left(|x_1-x_2|+\EE|Y^{x_1,0,m}_t-Y^{x_2,0,m}_t|\right)\\
\leq\!\!\!\!\!\!\!\!&&Ce^{-\frac{(\lambda_1-L_F) t}{2}}+C|x_1-x_2|.
\end{eqnarray*}
As a result, the proof is completed by letting $t\rightarrow \infty$.
\end{proof}
\begin{remark}\label{R4.2}
By a similar argument above, it is easy to prove that the averaged coefficient $\bar{B}$ is also Lipschitz continuous. As a consequence, equation \eref{1.3} admits a unique mild solution $\bar{X}_{t}$.
\end{remark}

Now we are in a position to prove Theorem \ref{main result 1}.

\vspace{0.1cm}
\noindent
\textbf{Proof of Theorem \ref{main result 1}}. It is easy to see that for any $T>0$, $p\in [1,\alpha)$ and $m\in \mathbb{N}_{+}$,
\begin{eqnarray*}
\sup_{t\in [0, T]}\EE|X_{t}^{\vare}-\bar{X}_{t}|^p\leq\!\!\!\!\!\!\!\!&&C_p\sup_{t\in [0, T]}\EE|X_{t}^{m,\vare}-X_{t}^{\vare}|^p+C_p\sup_{t\in [0, T]}\EE|\bar{X}^m_{t}-\bar{X}_{t}|^p\\
&&+C_p\sup_{t\in [0, T]}\EE|X_{t}^{m,\vare}-\bar{X}^m_{t}|^p.
\end{eqnarray*}
By Lemmas \ref{GA1} and \ref{GA2}, it follows
\begin{eqnarray*}
\lim_{m\rightarrow\infty}\sup_{t\in [0, T]}\EE|X_{t}^{m,\vare}-X_{t}^{\vare}|^p=0,\quad \lim_{m\rightarrow\infty}\sup_{t\in [0, T]}\EE|\bar{X}^m_{t}-\bar{X}_{t}|^p=0.
\end{eqnarray*}
Thus it is sufficient to prove the for any $(x,y)\in H^{\eta}\times H$ with $\eta\in (0,1)$, $T>0$ and small enough $\vare,\delta>0$, there exists a positive constant $C_{p,T,\delta}$
independent of $m$ such that
\begin{eqnarray}
\sup_{t\in [0, T]}\EE|X_{t}^{m,\vare}-\bar{X}^m_{t}|^p\leq C_{p,T,\delta}(1+\|x\|^{(1+\delta)p}_{\eta}+|y|^{(1+\delta)p})\vare^{p\left(1-\frac{1}{\alpha}\right)},\label{M1}
\end{eqnarray}
which will proved by the following three steps.

\vspace{2mm}
\textbf{Step 1.}
Using the formulation of the mild solutions $X_{t}^{m,\vare}$ and $\bar{X}^{m}_{t}$, we have for any $t>0$,
\begin{eqnarray*}
X_{t}^{m,\vare}\!-\!\bar{X}^{m}_{t}=\!\!\!\!\!\!\!\!&&\int_{0}^{t}e^{(t-s)A}\left[B^m(X_{s}^{m,\vare},Y_{s}^{m,\vare})-\bar{B}^m(\bar{X}^{m}_{s})\right]ds\\
=\!\!\!\!\!\!\!\!&&\int_{0}^{t}\!\!\!e^{(t-s)A}\!\!\left[B^m( X^{m,\vare}_{s},Y_{s}^{m,\vare})\!-\!\bar{B}^m(X^{m,\vare}_{s})\right]ds+\!\int_{0}^{t}\!\!\!e^{(t-s)A}\!\!
\left[\bar{B}^m(X^{m,\vare}_{s})\!-\!\bar{B}^m(\bar{X}^m_{s})\right]ds.
\end{eqnarray*}

Note that the averaged coefficient $\bar{B}^m$ has been proved that it is Lipschitz continuous in Lemma \ref{barX}. Then we get for any $T>0$,
\begin{eqnarray*}
\sup_{t\in [0, T]}\EE|X_{t}^{m,\vare}-\bar{X}^m_{t}|^p\leq\!\!\!\!\!\!\!\!&&C_p\sup_{t\in[0,T]}\EE\left|\int_{0}^{t}e^{(t-s)A}\left[B^m(X^{m,\vare}_{s},Y_{s}^{m,\vare})
-\bar{B}^m(X^{m,\vare}_{s})\right]ds\right|^p\\
&&+C_{p,T}\EE\int_{0}^{T}|X_{t}^{m,\vare}-\bar{X}^m_{t}|^p dt .
\end{eqnarray*}
By Gronwall's inequality, it follows
\begin{eqnarray}
\sup_{t\in [0, T]}\EE|X_{t}^{m,\vare}-\bar{X}^m_{t}|^p\leq C_{p,T}\sup_{t\in[0,T]}\EE\left|\int_{0}^{t}e^{(t-s)A}\left[B^m(X^{m,\vare}_{s},Y_{s}^{m,\vare})
-\bar{B}^m(X^{m,\vare}_{s})\right]ds\right|^p.\label{ID1}
\end{eqnarray}

By Proposition \ref{P3.6}, the following Poisson equation
\begin{eqnarray}
-\mathscr{L}^m_{2}(x)\Phi_m(x,y)=B^m(x,y)-\bar{B}^m(x)\label{PE}
\end{eqnarray}
exists a solution $\Phi_m(x,y)$ satisfies \eref{E1}-\eref{E3}.
\vspace{2mm}

By applying It\^o's formula,
\begin{eqnarray}
\Phi_m(X^{m,\vare}_{t},Y^{m,\vare}_{t})=\!\!\!\!\!\!\!\!&&e^{tA}\Phi_m(x^m,y^m)+\int^t_0 (-A)e^{(t-s)A}\Phi_m(X^{m,\vare}_{s},Y^{m,\vare}_{s})ds\nonumber\\
&&+\int^t_0 e^{(t-s)A}\mathscr{L}^m_{1}(Y^{m,\vare}_{s})\Phi_m(X^{m,\vare}_{s},Y^{m,\vare}_{s})ds+M^{m,\vare,1}_{t}+M^{m,\vare,2}_{t}\nonumber\\
&&+\frac{1}{\vare}\int^t_0 e^{(t-s)A}\mathscr{L}^m_{2}(X_{s}^{m,\vare})\Phi_m(X^{m,\vare}_{s},Y^{m,\vare}_{s})ds,\label{ID3}
\end{eqnarray}
where $\mathscr{L}^m_{1}(y)$, $M^{m,\vare,1}_{t}$ and $M^{m,\vare,2}_{t}$ are defined as follows:
\begin{eqnarray*}
\mathscr{L}^m_{1}(y)\Phi_m(x,y):=\!\!\!\!\!\!\!\!&&D_x\Phi_m(x,y)\cdot(Ax+B^m(x,y))\\
&&+\sum^m_{k=1}\beta^{\alpha}_k \int_{\RR}\Phi_m(x+e_k z,y)-\Phi_m(x,y)- D_x\Phi_m(x,y)\cdot (e_k z) 1_{\{|z|\leq 1\}}\nu (dz)
\end{eqnarray*}
$$M^{m,\vare,1}_{t}:=\sum^m_{k= 1}\int^t_0 \int_{\RR} e^{(t-s)A}\left[\Phi_m(X^{m,\vare}_{s-}+z\beta_k e_k, Y^{\vare}_{s-})-\Phi_m(X^{m,\vare}_{s-}, Y^{m,\vare}_{s-})\right]\tilde{N}^{1,k}(ds,dz),$$
$$M^{m,\vare,2}_{t}:=\sum^m_{k= 1}\int^t_0 \int_{\RR}e^{(t-s)A}\left[\Phi_m(X^{m,\vare}_{s-}, Y^{m,\vare}_{s-}+\vare^{-1/\alpha}z\gamma_k e_k)
-\Phi_m(X^{m,\vare}_{s-}, Y^{m,\vare}_{s-})\right]\tilde{N}^{2,k}(ds,dz).$$

Then it follows
\begin{eqnarray}
&&\int^t_0 -e^{(t-s)A}\mathscr{L}^m_{2}(X_{s}^{m,\vare})\Phi_m(X^{m,\vare}_{s},Y^{m,\vare}_{s}) ds\nonumber\\
=\!\!\!\!\!\!\!\!&&\vare\Big[e^{tA}\Phi_m(x^m,y^m)-\Phi_m(X_{t}^{m,\vare},Y^{m,\vare}_{t})+\int^t_0 (-A)e^{(t-s)A}\Phi_m(X^{m,\vare}_{s},Y^{m,\vare}_{s})ds\nonumber\\
&&+M^{m,\vare,1}_{t}+M^{m,\vare,2}_{t}+\int^t_0 \!\!e^{(t-s)A}\mathscr{L}^m_{1}(Y^{m,\vare}_{s})\Phi_m(X^{m,\vare}_{s},Y^{m,\vare}_{s})ds\Big].\label{ID4}
\end{eqnarray}
By \eref{ID1}, \eref{PE} and \eref{ID4}, we get
\begin{eqnarray}
\sup_{t\in [0, T]}\EE|X_{t}^{m,\vare}-\bar{X}^m_{t}|^p
\leq\!\!\!\!\!\!\!\!&&C_{p,T}\vare^{p}\Bigg\{\sup_{t\in[0,T]}\EE|e^{tA}\Phi_m(x^m,y^m)-\Phi_m(X_{t}^{m,\vare},Y^{m,\vare}_{t})|^p\nonumber\\
&&+\sup_{t\in[0,T]}\EE\left|\int^t_0 (-A)e^{(t-s)A}\Phi_m(X^{m,\vare}_{s},Y^{m,\vare}_{s})ds\right|^p\nonumber\\
&&+\sup_{t\in[0,T]}\EE|M^{m,\vare,1}_{t}|^p+\sup_{t\in[0,T]}\EE|M^{m,\vare,2}_{t}|^p\nonumber\\
&&+\sup_{t\in[0,T]}\EE\left|\int^t_0 e^{(t-s)A}\mathscr{L}^m_{1}(Y^{m,\vare}_{s})\Phi(X_{s}^{m,\vare},Y^{m,\vare}_{s})ds\right|^p\Bigg\}\nonumber\\
:=\!\!\!\!\!\!\!\!&&C_{p,T}\vare^{p}\sum^5_{k=1}\Lambda^{\vare,m}_k(T).\label{F5.4}
\end{eqnarray}

\textbf{Step 2.} In this step, we first estimate the terms $\Lambda^{\vare,m}_1(T)$-$\Lambda^{\vare,m}_4(T)$.
By \eref{E1}, \eref{Xvare} and \eref{Yvare}, we have
\begin{eqnarray}
\Lambda^{\vare,m}_1(T)\leq\!\!\!\!\!\!\!\!&&C\left(1+|x|^p+|y|^p+\sup_{t\in[0,T]}\EE|X^{m,\vare}_{t}|^p+\sup_{t\in[0,T]}\EE|Y^{m,\vare}_{t}|^p\right)\nonumber\\
\leq\!\!\!\!\!\!\!\!&&C_T(1+|x|^p+|y|^p).\label{F5.5}
\end{eqnarray}

By \eref{E1}, \eref{Xvare}, \eref{Yvare} and Lemma \ref{L6.3}, for any $\eta\in (0,1)$ we have
\begin{eqnarray}
\Lambda^{\vare,m}_2(T)\leq\!\!\!\!\!\!\!\!&&C_p\sup_{t\in[0,T]}\EE\left|\int^t_0 (-A)e^{(t-s)A}\left[\Phi_m(X^{m,\vare}_{s},Y^{m,\vare}_{s})
-\Phi_m(X^{m,\vare}_{t},Y^{m,\vare}_{t})\right]ds\right|^p\nonumber\\
&&+C_p\sup_{t\in[0,T]}\EE\left|\int^t_0 (-A)e^{(t-s)A}\Phi_m(X^{m,\vare}_{t},Y^{m,\vare}_{t})ds\right|^p\nonumber\\
\leq\!\!\!\!\!\!\!\!&&C_p\sup_{t\in[0,T]}\left|\int^t_0 (t-s)^{-1}\left[\EE|\Phi_m(X^{m,\vare}_{s},Y^{m,\vare}_{s})-\Phi_m(X^{m,\vare}_{t},Y^{m,\vare}_{t})|^p\right]^{1/p}ds\right|^p\nonumber\\
&&+C_p\sup_{t\in[0,T]}\EE\left|\Phi_m(X^{m,\vare}_{t},Y^{m,\vare}_{t})\right|^p\nonumber\\
\leq\!\!\!\!\!\!\!\!&&C_p\sup_{t\in[0,T]}\left|\int^t_0 (t-s)^{-1}\left[\EE\left((1+|X^{m,\vare}_{s}|^{\delta p}+|X^{m,\vare}_{t}|^{\delta p}
+|Y^{m,\vare}_{s}|^{\delta p})|X^{m,\vare}_{s}-X^{m,\vare}_{t}|^p\right)\right]^{1/p}ds\right|^p\nonumber\\
&&+C_p\sup_{t\in[0,T]}\left|\int^t_0 (t-s)^{-1}\left[\EE|Y^{m,\vare}_{s}-Y^{m,\vare}_{t}|^p\right]^{1/p}ds\right|^p\nonumber\\
&&+C_p\sup_{t\in[0,T]}\EE\left|\Phi_m(X^{m,\vare}_{t},Y^{m,\vare}_{t})\right|^p\nonumber\\
\leq\!\!\!\!\!\!\!\!&&C_{p,T}(1+|x|^{p(1+\delta)}+|y|^{p(1+\delta)})\sup_{t\in[0,T]}\left|\int^t_0 (t-s)^{-1}(t-s)^{\frac{\eta}{2}}s^{-\frac{\eta}{2}}ds\right|^p\nonumber\\
&&+C_{p,T}(1+|x|^p+|y|^{p})\sup_{t\in[0,T]}\left|\int^t_0 (t-s)^{-1}\left(\frac{t-s}{\vare}\right)^{\frac{\eta}{2}}s^{-\frac{\eta}{2}}ds\right|^p\nonumber\\
&&+C_{p,T}(1+|x|^p+|y|^{p})\nonumber\\
\leq\!\!\!\!\!\!\!\!&& C_{p,T}(1+|x|^{p(1+\delta)}+|y|^{p(1+\delta)})\vare^{-\frac{\eta p}{2}}.\label{F5.6}
\end{eqnarray}

For the term $\Lambda^{\vare,m}_3(T)$. By Burkholder-Davis-Gundy's inequality, it follows for any $p\in [1,\alpha)$,

\begin{eqnarray}
\Lambda^{\vare,m}_3(T)
\leq\!\!\!\!\!\!\!\!&&C_p\EE\left[\sum^m_{k=1}\int^T_0 \int_{\RR}|\Phi_m(X^{m,\vare}_{s-}\!+\!z\beta_k e_k, Y^{m,\vare}_{s-})\!
-\!\Phi_m(X^{m,\vare}_{s-}, Y^{m,\vare}_{s-})|^2 N^{1,k}(ds,dz)\!\right]^{p/2}\nonumber\\
\leq\!\!\!\!\!\!\!\!&&C_p\EE\left[\sum^m_{k=1}\int^T_0 \int_{|z|\leq \beta^{-1}_k}|\Phi_m(X^{m,\vare}_{s-}\!+\!z\beta_k e_k, Y^{m,\vare}_{s-})\!
-\!\Phi_m(X^{m,\vare}_{s-}, Y^{m,\vare}_{s-})|^2 N^{1,k}(ds,dz)\!\right]^{p/2}\nonumber\\
&&+C_p\EE\left[\sum^m_{k=1}\int^T_0 \int_{|z|>\beta^{-1}_k}|\Phi_m(X^{m,\vare}_{s-}+z\beta_k e_k, Y^{m,\vare}_{s-})-\Phi_m(X^{m,\vare}_{s-}, Y^{m,\vare}_{s-})|^2 N^{1,k}(ds,dz)\right]^{p/2}\nonumber\\
\leq\!\!\!\!\!\!\!\!&&C_p\left[\EE\sum^m_{k=1}\int^T_0 \int_{|z|\leq \beta^{-1}_k}|\Phi_m(X^{m,\vare}_{s-}\!+\!z\beta_k e_k, Y^{m,\vare}_{s-})\!
-\!\Phi_m(X^{m,\vare}_{s-}, Y^{m,\vare}_{s-})|^2 N^{1,k}(ds,dz)\!\right]^{p/2}\nonumber\\
&&+C_p\EE\sum^m_{k=1}\int^T_0 \int_{|z|>\beta^{-1}_k}|\Phi_m(X^{m,\vare}_{s-}+z\beta_k e_k, Y^{m,\vare}_{s-})-\Phi_m(X^{m,\vare}_{s-}, Y^{m,\vare}_{s-})|^p N^{1,k}(ds,dz)\nonumber\\
\leq\!\!\!\!\!\!\!\!&&C_p\left[\sum^m_{k=1}\beta^{\alpha}_k\EE\int^T_0 \int_{|z|\leq 1}|\Phi_m(X^{m,\vare}_{s}+z e_k, Y^{m,\vare}_{s})
-\Phi_m(X^{m,\vare}_{s}, Y^{m,\vare}_{s})|^2 \nu(dz)ds\right]^{p/2}\nonumber\\
&&+C_p\sum^m_{k=1}\beta^{\alpha}_k\EE\int^T_0 \int_{|z|> 1}|\Phi_m(X^{m,\vare}_{s}+z e_k, Y^{m,\vare}_{s})-\Phi_m(X^{m,\vare}_{s}, Y^{m,\vare}_{s})|^p \nu(dz)ds.\nonumber
\end{eqnarray}
By \eref{E2}, for any $\delta\in (0,1]$ we have
\begin{eqnarray*}
|\Phi_m(X^{m,\vare}_{s}+z e_k, Y^{m,\vare}_{s})-\Phi_m(X^{m,\vare}_{s}, Y^{m,\vare}_{s})|=\!\!\!\!\!\!\!\!&&\left|\int^1_0 D_x\Phi_m(X^{m,\vare}_{s}+\xi z e_k, Y^{m,\vare}_{s})\cdot (z e_k) d\xi \right|\\
\leq\!\!\!\!\!\!\!\!&&C_{\delta}\int^1_0(1+|X^{m,\vare}_{s}+\xi z e_k|^{\delta}+|Y^{m,\vare}_{s}|^{\delta})d \xi |z e_k|\nonumber\\
\leq\!\!\!\!\!\!\!\!&&C_{\delta}(1+|X^{m,\vare}_{s}|^{\delta}+|Y^{m,\vare}_{s}|^{\delta}+|z|^{\delta})|z|.
\end{eqnarray*}
Note that for small enough $\delta<\frac{\alpha}{p}-1$, it easy to see
$$
\int_{|z|\leq 1}|z|^2\nu(dz)<\infty,\quad \int_{|z|>1}|z|^{p(1+\delta)} \nu(dz)<\infty,
$$
which combing with the condition $\sum^{\infty}_{k=1} \beta^{\alpha}_k<\infty$ we obtain
\begin{eqnarray}
\Lambda^{\vare,m}_3(T)
\leq\!\!\!\!\!\!\!\!&&C_{p,\delta}\left[\sum^m_{k=1}\beta^{\alpha}_k\int^T_0 \int_{|z|\leq 1}\EE(1+|X^{m,\vare}_{s}|^{2\delta}+|Y^{m,\vare}_{s}|^{2\delta}+|z|^{2\delta})|z|^2 \nu(dz)ds\right]^{p/2}\nonumber\\
&&+C_{p,\delta}\sum^m_{k=1}\beta^{\alpha}_k\EE\int^T_0 \int_{|z|> 1}\EE(1+|X^{m,\vare}_{s}|^{p\delta}+|Y^{m,\vare}_{s}|^{p\delta}+|z|^{p\delta})|z|^p \nu(dz)ds\nonumber\\
\leq\!\!\!\!\!\!\!\!&&C_{p,T,\delta}(1+|x|^{p\delta}+|y|^{p\delta}).\label{F5.7}
\end{eqnarray}

For the term $\Lambda^{\vare,m}_4(T)$. Similar as we did in estimating $\Lambda^{\vare,m}_3(T)$. It is easy to see that
\begin{eqnarray}
\Lambda^{\vare,m}_4(T)\leq\!\!\!\!\!\!\!\!&&C_p\EE\left[\sum^m_{k=1}\int^T_0 \int_{\RR}|\Phi_m(X^{m,\vare}_{s-},Y^{m,\vare}_{s-}\!+\!\vare^{-1/\alpha}\gamma_k z e_k)\!
-\!\Phi_m(X^{m,\vare}_{s-},Y^{m,\vare}_{s-})|^2 N^{2,k}(ds,dz)\!\right]^{p/2}\nonumber\\
\leq\!\!\!\!\!\!\!\!&&C_p\EE\left[\sum^m_{k=1}\!\!\int^T_0 \!\!\int_{|z|\leq \gamma^{-1}_k}|\Phi_m(X^{m,\vare}_{s-},Y^{m,\vare}_{s-}\!+\!\vare^{-1/\alpha}\gamma_k z e_k)\!
-\!\Phi_m(X^{m,\vare}_{s-},Y^{m,\vare}_{s-})|^2 N^{2,k}(ds,dz)\!\right]^{p/2}\nonumber\\
&&\!\!\!\!\!\!\!\!+C_p\EE\left[\sum^m_{k=1}\!\!\int^T_0 \int_{|z|>\gamma^{-1}_k}|\Phi_m(X^{m,\vare}_{s-},Y^{m,\vare}_{s-}\!+\!\vare^{-1/\alpha}\gamma_k z e_k)\!
-\!\Phi_m(X^{m,\vare}_{s-},Y^{m,\vare}_{s-})|^2 N^{2,k}(ds,dz)\right]^{p/2}\nonumber\\
\leq\!\!\!\!\!\!\!\!&&C_p\left[\EE\sum^m_{k=1}\!\!\int^T_0\!\int_{|z|\leq \gamma^{-1}_k}|\Phi_m(X^{m,\vare}_{s-},Y^{m,\vare}_{s-}\!+\!\vare^{-1/\alpha}\gamma_k z e_k)\!
-\!\Phi_m(X^{m,\vare}_{s-},Y^{m,\vare}_{s-})|^2 N^{2,k}(ds,dz)\!\right]^{p/2}\nonumber\\
&&\!\!\!\!\!\!\!\!+C_p\EE\sum^m_{k=1}\!\!\int^T_0 \int_{|z|>\gamma^{-1}_k}|\Phi_m(X^{m,\vare}_{s-},Y^{m,\vare}_{s-}\!+\!\vare^{-1/\alpha}\gamma_k z e_k)\!
-\!\Phi_m(X^{m,\vare}_{s-},Y^{m,\vare}_{s-})|^p N^{2,k}(ds,dz)\nonumber\\
\leq\!\!\!\!\!\!\!\!&&C_p\left[\sum^m_{k=1}\gamma^{\alpha}_k\EE\int^T_0\! \int_{|z|\leq 1}|\Phi_m(X^{m,\vare}_{s},Y^{m,\vare}_{s}+\vare^{-1/\alpha}z e_k)
-\Phi_m(X^{m,\vare}_{s},Y^{m,\vare}_{s})|^2 \nu(dz)ds\right]^{p/2}\nonumber\\
&&\!\!\!\!\!\!\!\!+C_p\sum^m_{k=1}\gamma^{\alpha}_k\EE\int^T_0 \int_{|z|> 1}|\Phi_m(X^{m,\vare}_{s},Y^{m,\vare}_{s}+\vare^{-1/\alpha}z e_k)-\Phi_m(X^{m,\vare}_{s},Y^{m,\vare}_{s})|^p \nu(dz)ds.\nonumber
\end{eqnarray}
By \eref{E21}, we have
\begin{eqnarray*}
&&|\Phi_m(X^{m,\vare}_{s},Y^{m,\vare}_{s}+\vare^{-1/\alpha}z e_k)-\Phi_m(X^{m,\vare}_{s},Y^{m,\vare}_{s})|\\
=\!\!\!\!\!\!\!\!&&\left|\int^1_0 D_y\Phi_m(X^{m,\vare}_{s}, Y^{m,\vare}_{s}+\xi\vare^{-1/\alpha}z e_k)\cdot (\vare^{-1/\alpha}z e_k) d\xi \right|\leq C\vare^{-1/\alpha}|z|.
\end{eqnarray*}
Then by condition $\sum^{\infty}_{k=1} \gamma^{\alpha}_k<\infty$ we obtain
\begin{eqnarray}
\Lambda^{\vare,m}_4(T)
\leq\!\!\!\!\!\!\!\!&&C_p\vare^{-p/\alpha}\left[\sum^m_{k=1}\gamma^{\alpha}_k\int^T_0 \int_{|z|\leq 1}|z|^2 \nu(dz)ds\right]^{p/2}\nonumber\\
&&+C_p\vare^{-p/\alpha}\sum^m_{k=1}\gamma^{\alpha}_k\EE\int^T_0 \int_{|z|> 1}|z|^p \nu(dz)ds\nonumber\\
\leq\!\!\!\!\!\!\!\!&&C_{p,T}\vare^{-p/\alpha}.\label{F5.8}
\end{eqnarray}

\textbf{Step 3.} In this step, we estimate the term $\Lambda^{\vare,m}_5(T)$. It is easy to see
\begin{eqnarray*}
\Lambda^{\vare,m}_5(T)\leq\!\!\!\!\!\!\!\!&&C_p\sup_{t\in[0,T]}\EE\left|\int^t_0 e^{(t-s)A}\left[D_x\Phi_m(X_{s}^{m,\vare},Y^{m,\vare}_{s})\cdot AX_{s}^{m,\vare}\right]ds\right|^p\nonumber\\
&&+C_p\sup_{t\in[0,T]}\EE\left|\int^t_0 e^{(t-s)A}\left[D_x\Phi_m(X_{s}^{m,\vare},Y^{m,\vare}_{s})\cdot B^m(X_{s}^{m,\vare},Y^{m,\vare}_{s})\right] ds\right|^p\nonumber\\
&&+C_p\sup_{t\in[0,T]}\EE\left|\sum^m_{k=1}\beta^{\alpha}_k \int^t_0\int_{\RR}e^{(t-s)A}\big[\Phi_m(x+e_k z,y)-\Phi_m(x,y)\right.\\
&&\quad\quad\quad\quad\quad\quad\quad\left.-D_x\Phi_m(x,y)\cdot( e_k z)\big]\nu (dz)ds\right|^p:=\sum^3_{i=1}\Lambda_{5i}(T).
\end{eqnarray*}

By \eref{Xvare2}, \eref{E2}, \eref{LA}, \eref{Xvare}, \eref{Yvare} and Minkowski's inequality, we get for any $\gamma\in (0,1)$ and $\eta\in (0,1)$,
\begin{eqnarray*}
\Lambda^{\vare,m}_{51}(T)\leq\!\!\!\!\!\!\!\!&&C_{p,\delta}\sup_{t\in[0,T]}\EE\left|\int^t_0 \|X_{s}^{m,\vare}\|_{2}(1+|X^{m,\vare}_{s}|^{\delta}+|Y^{m,\vare}_{s}|^{\delta} ) ds\right|^p\nonumber\\
\leq\!\!\!\!\!\!\!\!&&C_{p,\delta}\sup_{t\in[0,T]}\left|\int^t_0 \left[\EE\left(\|X^{m,\varepsilon}_s\|^p_{2}(1+|X^{m,\vare}_{s}|^{\delta p}
+|Y^{m,\vare}_{s}|^{\delta p} ) \right)\right]^{1/p}ds\right|^p\nonumber\\
\leq\!\!\!\!\!\!\!\!&&C_{p,\delta}\sup_{t\in[0,T]}\left|\int^t_0 \left[\EE\|X^{m,\varepsilon}_s\|^{p'}_{2}\right]^{1/p'}\left[\EE(1+|X^{m,\vare}_{s}|^{\frac{\delta pp'}{p'-p}}
+|Y^{m,\vare}_{s}|^{\frac{\delta pp'}{p'-p}} )\right]^{\frac{p'-p}{pp'}}ds\right|^p\nonumber\\
\leq\!\!\!\!\!\!\!\!&&C_{p,T,\delta}(1+\|x\|^{(1+\delta)p}_{\gamma}+|y|^{(1+\delta)p})\sup_{t\in[0,T]}\left|\int^t_0 s^{-1+\frac{\gamma}{2}}ds+\vare^{-\frac{\eta}{2}}\right|^p\nonumber\\
\leq\!\!\!\!\!\!\!\!&&C_{p,T,\delta}(1+\|x\|^{(1+\delta)p}_{\gamma}+|y|^{(1+\delta)p})\vare^{-\frac{\eta p}{2}},
\end{eqnarray*}
where $p<p'<\alpha$ and $\delta$ is small enough such that $\frac{\delta p'}{p'-p}\leq 1$.

By a similar argument above, we have for $\delta$ is small enough such that $\delta<\frac{\alpha}{p}-1$,
\begin{eqnarray*}
\Lambda^{\vare,m}_{52}(T)\leq\!\!\!\!\!\!\!\!&&C_{p,\delta}\EE\int^T_0 |B^m(X_{s}^{m,\vare},Y^{m,\vare}_{s})|^p (1+|X^{m,\vare}_{s}|^{\delta p}+|Y^{m,\vare}_{s}|^{\delta p})  ds\nonumber\\
\leq\!\!\!\!\!\!\!\!&&C_{p,\delta}\EE\int^T_0(1+|X_{s}^{m,\vare}|^p+|Y^{m,\vare}_{s}|^p)(1+|X^{m,\vare}_{s}|^{\delta p}+|Y^{m,\vare}_{s}|^{\delta p}) ds\\
\leq\!\!\!\!\!\!\!\!&&C_{p,T,\delta}(1+|x|^{(1+\delta)p}+|y|^{(1+\delta)p}).
\end{eqnarray*}

By \eref{E2} and \eref{E3}, we have for $\delta$ is small enough such that $\delta<\alpha-1$,
\begin{eqnarray*}
\Lambda^{\vare,m}_{53}(T)\leq\!\!\!\!\!\!\!\!&&C_p\EE\Big[\sum^m_{k=1}\beta^{\alpha}_k\int^T_0\int_{|z|\leq \tilde{c}_k}
|\Phi_m(X^{m,\vare}_{s}+e_k z,Y^{m,\vare}_{s})-\Phi_m(X^{m,\vare}_{s},Y^{m,\vare}_{s})\\
&&\quad\quad\quad\quad-D_x \Phi_m(X^{m,\vare}_{s},Y^{m,\vare}_{s})\cdot(e_k z)|\nu(dz)ds\Big]^p\nonumber\\
&&+C_{p}\EE\left[\sum^m_{k=1}\beta^{\alpha}_k\int^T_0\int_{|z|> \tilde{c}_k}|\Phi_m(X^{m,\vare}_{s}+e_k z,Y^{m,\vare}_{s})-\Phi_m(X^{m,\vare}_{s},Y^{m,\vare}_{s})|\nu(dz)ds\right]^p\nonumber\\
\leq\!\!\!\!\!\!\!\!&&C_{p}\EE\left[\sum^{\infty}_{k=1}\beta^{\alpha}_k\lambda^{\frac{\kappa_1}{2}}_k\int^T_0\int_{|z|\leq  \tilde{c}_k}|z|^2\nu(dz)(1+|X^{m,\vare}_{s}|+|Y^{m,\vare}_s|)ds\right]^p \nonumber\\
&&+C_{p,\delta}\EE\left[\sum^{\infty}_{k=1}\beta^{\alpha}_k\int^T_0\int_{|z|> \tilde{c}_k}|z|^{1+\delta}(1+|X^{m,\vare}_{s}|^{\delta}+|Y^{m,\vare}_s|^{\delta})\nu(dz)ds\right]^p\nonumber\\
\leq\!\!\!\!\!\!\!\!&&C_{p}\left(\sum^{\infty}_{k=1}\beta^{\alpha}_k\lambda^{\frac{\kappa_1}{2}}_k\tilde{c}^{2-\alpha}_k\right)^p\EE\left|\int^T_0(1+|X^{m,\vare}_{s}|+|Y^{m,\vare}_s|)ds\right|^p \nonumber\\
&&+C_{p,\delta}\left(\sum^{\infty}_{k=1}\beta^{\alpha}_k\tilde{c}^{\delta-\alpha+1}_k\right)^p\EE\left|\int^T_0(1+|X^{m,\vare}_{s}|^{\delta}+|Y^{m,\vare}_s|^{\delta})\nu(dz)ds\right|^p,
\end{eqnarray*}
where $\tilde{c}_k:=\lambda^{-\frac{\kappa_1}{2(1-\delta)}}_k$. Then by $\sum_{k\in \mathbb{N}_{+}}\frac{\beta^{\alpha}_k}{\lambda^{1-\alpha}_k}<\infty$, it follows
\begin{eqnarray*}
\Lambda^{\vare,m}_{53}(T)\leq\!\!\!\!\!\!\!\!&&C_{p,T,\delta}(1+|x|+|y|)\left(\sum^{\infty}_{k=1}\beta^{\alpha}_k\lambda^{\frac{\kappa_1 (\alpha-\delta-1)}{2(1-\delta)}}_k\right)^p\nonumber\\
\leq\!\!\!\!\!\!\!\!&&C_{p,T,\delta}(1+|x|^p+|y|^p)\left(\sum^{\infty}_{k=1}\beta^{\alpha}_k\lambda^{\frac{\alpha-\delta-1}{(1-\delta)}}_k\right)^p\nonumber\\
\leq\!\!\!\!\!\!\!\!&&C_{p,T,\delta}(1+|x|^p+|y|^p)\left(\sum^{\infty}_{k=1}\beta^{\alpha}_k\lambda^{\alpha-1}_k\right)^p\nonumber\\
\leq\!\!\!\!\!\!\!\!&&C_{p,T,\delta}(1+|x|^p+|y|^p).\label{F5.6}
\end{eqnarray*}

As a result, for small enough $\delta>0$ we have
\begin{eqnarray}
\Lambda^{\vare,m}_{5}(T)\leq C_{p,T,\delta}(1+|x|^{(1+\delta)p}+|y|^{(1+\delta)p}). \label{Gamma3}
\end{eqnarray}

Hence, \eref{M1} holds by combining \eref{F5.5}-\eref{Gamma3}. The proof is complete.

\section{Weak convergence order}

This section is devoted to prove Theorem \ref{main result 2}. We still consider the problem in finite dimension firstly, then passing the limit to the infinite dimensional case.

For a test function $\phi \in C_{b}^{3}(H)$, we have for any $t\geq 0$,
\begin{eqnarray} \label{ephix}
\left|\mathbb{E}\phi\left(X^{\vare}_t\right)
-\EE\phi(\bar{X}_t)\right|
\leq\!\!\!\!\!\!\!\!&&
\left|\mathbb{E}\phi\left(X^{\vare}_t\right)
-\mathbb{E}\phi\left(X^{m,\vare}_t\right)\right|+\left|\EE\phi (\bar{X}^m_t)-\EE\phi(\bar{X}_t)\right| \nonumber\\
\!\!\!\!\!\!\!\!&&+\left|\mathbb{E}\phi\left(X^{m,\vare}_t\right)-\EE\phi (\bar{X}^m_t)\right|.
\end{eqnarray}
By Lemmas \ref{GA1} and \ref{GA2} in the appendix, it is easy to see that
\begin{eqnarray*}
&&\lim_{m\rightarrow\infty}\sup_{t\in [0, T]}\left|\mathbb{E}\phi\left(X^{\vare}_t\right)
-\mathbb{E}\phi\left(X^{m,\vare}_t\right)\right|=0,\\
&&\lim_{m\rightarrow\infty}\sup_{t\in [0, T]}\left|\EE\phi (\bar{X}^m_t)-\EE\phi(\bar{X}_t)\right|=0.
\end{eqnarray*}
Then the proof will be completed if we can show that there exists a positive constant $C$ independent of $m$ such that
$$
\sup_{t\in [0, T]}\left|\mathbb{E}\phi\left(X^{m,\vare}_t\right)-\EE\phi (\bar{X}^m_t)\right|\leq C\vare.
$$

Now we are in a position to prove Theorem \ref{main result 2}.

\vspace{0.2cm}

\noindent
Proof of Theorem \ref{main result 2}. We will divided the proof into three steps.

\vspace{2mm}
\textbf{Step 1.} We first introduce the following Kolmogorov equation in finite dimension:
\begin{equation}\left\{\begin{array}{l}\label{KE}
\displaystyle
\partial_t u_m(t,x)=\bar{\mathscr{L}}^m_1 u_m(t,x),\quad t\in[0, T], \vspace{2mm}\\
u(0, x)=\phi(x),\quad x\in H_m,
\end{array}\right.
\end{equation}
where $\phi\in C^{3}_b(H)$ and $\bar{\mathscr{L}}^m_1$ is the infinitesimal generator of the transition semigroup of the averaged equation \eref{Ga 1.3}, which is given by
\begin{eqnarray*}
\bar{\mathscr{L}}^m_{1}\phi(x):=\!\!\!\!\!\!\!\!&& D_x\phi(x)\cdot [Ax+\bar{B}^m(x)]\nonumber\\
&&+\sum^{m}_{k=1}\beta^{\alpha}_k \int_{\RR}\left[\phi(x+e_k z)-\phi(x)-\langle D_x\phi(x), e_k z\rangle 1_{\{|z|\leq 1\}}\right]\nu (dz).
\end{eqnarray*}

Note that
$$\bar{B}^m(x):=\int_{H_m}B^m(x,y)\mu^{x,m}(dy)=\lim_{t\rightarrow \infty}\EE B^m(x, Y^{x,y,m}_t).$$
By assumptions \ref{A3}, conditions \eref{ThirdDer F4} and \eref{ThirdDer F5}, then through a straightforward computation, it is easy to check that
$$|D\bar{B}^m(x)\cdot h|\leq C |h| \quad \forall x, h\in H_m ,$$
$$|D^2\bar{B}^m(x)\cdot(h,k)|\leq C |h|\|k\|_{\kappa_1},\quad \forall x, h, k\in H_m,$$
$$|D^3\bar{B}^m(x)\cdot(h,k,l)|\leq C|h|\|k\|_{\kappa_1}\|l\|_{\kappa_1}, \quad \forall x,h,k,l\in H_m,$$
where $\kappa_1$ is the constant in assumption \ref{A3}. As a consequence, \eref{KE} has a unique solution $u_m$ given by
$$
u_m(t,x)=\EE\phi(\bar{X}^m_t(x)),\quad t\in [0,T].
$$
Furthermore, for any $h,k,l\in H_m$ we have
\begin{eqnarray}
\sup_{x\in H_m,m\geq 1}|D_{x} u_m(t,x)\cdot h|\leq C_T |h|,\quad t\in [0,T],\label{UE3}
\end{eqnarray}
\begin{eqnarray}
\sup_{x\in H_m,m\geq 1}|D_{xx} u_m(t,x)\cdot (h,k)|\leq C_T |h||k|,\quad t\in[0,T].\label{UE4}
\end{eqnarray}
\begin{eqnarray}
\sup_{x\in H_m,m\geq 1}|D_{xxx} u_m(t,x)\cdot (h,k,l)|\leq C_T|h||k||l|,\quad t\in[0,T],\label{UE1}
\end{eqnarray}
\begin{eqnarray}
\sup_{m\geq 1}|\partial_t(D_x u_m(t,x))\cdot h|\leq C_T t^{-1}(1+|x|)|h|,\quad x\in H_m, t\in(0,T].\label{UE2}
\end{eqnarray}

Indeed, note that for any $h,k,l\in H_m$,
$$D_x u_m(t,x)\cdot h=\EE[D\phi(\bar{X}^m_t)\cdot \eta^{h,m}_t(x)],$$
$$D_{xx} u_m(t,x)\cdot (h,k)=\EE\left[D^2\phi(\bar{X}^m_t)\cdot (\eta^{h,m}_t(x),\eta^{k,m}_t(x))\right]+\EE\left[D\phi(\bar{X}^m_t)\cdot  \zeta^{h,k,m}_t(x)\right],$$
\begin{eqnarray*}
D_{xxx} u_m(t,x)\cdot (h,k,l)=\!\!\!\!\!\!\!\!&&\EE\left[D^3\phi(\bar{X}^m_t)\cdot (\eta^{h,m}_t(x),\eta^{k,m}_t(x),\eta^{l,m}_t(x))\right]\\
&&+\EE\left[D^2\phi(\bar{X}^m_t)\cdot (\zeta^{h,l,m}_t(x),\eta^{k,m}_t(x))+D^2\phi(\bar{X}^m_t)\cdot (\eta^{h,m}_t(x),\zeta^{k,l,m}_t(x))\right]\\
&&+\EE\left[D^2\phi(\bar{X}^m_t)\cdot  (\zeta^{h,k,m}_t(x), \eta^{l,m}_t(x))+D\phi(\bar{X}^m_t)\cdot  \chi^{h,k,l,m}_t(x)\right],
\end{eqnarray*}
where $\eta^{h,m}_t(x):=D_x \bar{X}^m_t(x)\cdot h$ satisfies
 \begin{equation}\left\{\begin{array}{l}
\displaystyle
d\eta^{h,m}_t(x)\!=\!A\eta^{h,m}_t(x)dt+D\bar{B}^m(\bar{X}^m_t)\cdot \eta^{h,m}_t(x)dt,\nonumber\\
\eta^{h,m}_0(x)=h,\nonumber
\end{array}\right.
\end{equation}
$\zeta^{h,k,m}_t(x):=D_{xx} \bar{X}^m_t(x)\cdot (h,k)$ satisfies
\begin{equation}\left\{\begin{array}{l}
\displaystyle
d\zeta^{h,k,m}_t(x)=\left[A\zeta^{h,k,m}_t(x)+\!\!D\bar{B}^m(\bar{X}^m_t)\cdot \zeta^{h,k,m}_t(x)+D^2\bar{B}^m(\bar{X}^m_t)\cdot (\eta^{h,m}_t(x), \eta^{k,m}_t(x))\right]\!\!dt,\nonumber\\
\zeta^{h,k,m}_0(x)=0\nonumber
\end{array}\right.
\end{equation}
and $\chi^{h,k,l,m}_t(x):=D_{xxx} \bar{X}^m_t(x)\cdot (h,k,l)$ satisfies
 \begin{equation}\left\{\begin{array}{l}
\displaystyle
d\chi^{h,k,l,m}_t(x)=\left[A\chi^{h,k,m}_t(x)+\!\!D\bar{B}^m(\bar{X}^m_t)\cdot \chi^{h,k,l,m}_t(x)+D^2\bar{B}^m(\bar{X}^m_t)\cdot (\zeta^{h,k,m}_t(x), \eta^{l,m}_t(x))\right.\nonumber\\
\quad\quad\quad\quad\quad\quad +D^2\bar{B}^m(\bar{X}^m_t)\cdot (\zeta^{h,l,m}_t(x),\eta^{k,m}_t(x))+D^2\bar{B}^m(\bar{X}^m_t)\cdot (\eta^{h,m}_t(x),\zeta^{k,l,m}_t(x))\nonumber\\
\quad\quad\quad\quad\quad\quad+\left.D^3\bar{B}^m(\bar{X}^m_t)\cdot (\eta^{h,m}_t(x), \eta^{k,m}_t(x), \eta^{l,m}_t(x))\right]dt,\nonumber\\
\chi^{h,k,l,m}_0(x)=0.\nonumber
\end{array}\right.
\end{equation}
By a straightforward computation, we obtain for any $t\in(0,T]$ and $\theta\in [0,2)$,
\begin{eqnarray}
&&\|\eta^{h,m}_t(x)\|_{\theta}\leq C_T t^{-\theta/2}|h|, \quad \|\eta^{h,m}_t(x)\|_{2}\leq C_T t^{-1}(1+|x|)|h|,\label{D barX}\\
&&|\zeta^{h,k,m}_t(x)|\leq C_T |h||k|,\label{Dxx barX}\\
&&|\chi^{h,k,l,m}_t(x)|\leq C_T|h||k||l|.\label{Dxxx barX}
\end{eqnarray}
Hence, it is easy to see \eref{D barX}-\eref{Dxxx barX} imply \eref{UE3}-\eref{UE1} hold.

By It\^{o}'s formula and taking expectation, we have
\begin{eqnarray*}
\EE[D\phi(\bar{X}^m_t)\cdot \eta^{h,m}_t(x)]=\!\!\!\!\!\!\!\!&&D\phi(x)\cdot h+\int^t_0 \EE\left[ D^2\phi(\bar{X}^m_s)\cdot (\eta^{h,m}_s(x), A\bar{X}^m_s+\bar{B}^m(\bar{X}^m_s))\right]ds\\
&&+\int^t_0 \EE \langle D\phi(\bar{X}^m_s), A\eta^{h,m}_s(x)+D\bar{B}^m(\bar{X}^m_s)\cdot \eta^{h,m}_s(x)\rangle ds\\
&&+\sum^m_{k=1}\beta^{\alpha}_k \EE\Big[\int^t_0\int_{\RR} \langle D\phi(\bar{X}^m_s+e_k z)-D\phi(\bar{X}^m_s)\\
&&\quad\quad\quad\quad-\left(D^2\phi(\bar{X}^m_s)\cdot e_k z\right) 1_{\{|z|\leq 1\}},  \eta^{h,m}_s(x)\rangle\nu(dz)ds\Big],
\end{eqnarray*}
which implies
\begin{eqnarray*}
\partial_t(D_x u_m(t,x))\cdot h=\!\!\!\!\!\!\!\!&& \EE\left[D^2\phi(\bar{X}^m_t)\cdot (\eta^{h,m}_t(x), A\bar{X}^m_t+\bar{B}^m(\bar{X}^m_t))\right]\\
&&+\EE \langle D\phi(\bar{X}^m_t), A\eta^{h,m}_t(x)+D\bar{B}^m(\bar{X}^m_t)\cdot \eta^{h,m}_t(x)\rangle\\
&&+\sum^m_{k=1}\beta^{\alpha}_k \EE\Big[\int_{\RR}\langle D\phi(\bar{X}^m_t+e_k z)-D\phi(\bar{X}^m_t)\\
&&\quad\quad\quad-\left[D^2\phi(\bar{X}^m_t)\cdot e_k z\right] 1_{\{|z|\leq 1\}},  \eta^{h,m}_t(x)\rangle\nu(dz)\Big].
\end{eqnarray*}
By \eref{Control X}, \eref{D barX} and \eref{X_theta}, we get for any $t\in (0,T]$,
\begin{eqnarray*}
|\partial_t(D_x u_m(t,x))\cdot h|\leq \!\!\!\!\!\!\!\!&& C\EE\left[\|\eta^{h,m}_t(x)\|_2+|\eta^{h,m}_t(x)|(\|\bar{X}^m_t\|_2+|\bar{X}^m_t|+1)\right]\\
&&+C\sum^m_{k=1}\beta^{\alpha}_k \EE|\eta^{h,m}_t(x)|\left[\int_{|z|\leq 1}|z|^2\nu(dz)+\int_{|z|>1}|z|\nu(dz)\right]\\
\leq\!\!\!\!\!\!\!\!&& C_T t^{-1}|h|(1+|x|),
\end{eqnarray*}
which proves \eref{UE2}.

\vspace{0.2cm}
\textbf{Step 2.} Let $\tilde{u}^t_m(s,x):=u_m(t-s,x)$, $s\in [0,t]$, by It\^{o}'s formula, we have
\begin{eqnarray*}
\tilde{u}^t_m(t, X^{m,\vare}_t)=\!\!\!\!\!\!\!\!&&\tilde{u}^t_m(0,x)+\int^t_0 \partial_s \tilde{u}^t_m(s, X^{m,\vare}_s )ds+\int^t_0 \mathscr{L}^m_{1}(Y^{m,\vare}_s)\tilde{u}^t_m(s, X^{m,\vare}_s)ds
+\tilde{M}^m_t,
\end{eqnarray*}
where $\tilde{M}^m_t$ is defined as follows,
$$
\tilde{M}^m_t:=\sum^{m}_{k=1}\beta^{\alpha}_k \int^t_0 \int_{\RR}\tilde{u}^t_m(s,X_{s-}^{m,\vare}+z e_k)-\tilde{u}^t_m(s,X_{s-}^{m,\vare})\tilde{N}^{1,k}(dz,ds).
$$
Note that
$$\tilde{u}^t_m(t, X^{m,\vare}_t)=\phi(X^{m,\vare}_t),\tilde{u}^t_m(0, x)=\EE\phi(\bar{X}^m_t(x)),\partial_s \tilde{u}^t_m(s, X^{m,\vare}_s )=-\bar{\mathscr{L}}^m_1 \tilde{u}^t_m(s, X^{m,\vare}_s).$$
It follows for any $t\in[0,T]$,
\begin{eqnarray*}
\left|\EE\phi(X^{m,\vare}_{t})-\EE\phi(\bar{X}^m_{t})\right|=\!\!\!\!\!\!\!\!&&\left|\EE\int^t_0 -\bar{\mathscr{L}}^m_1 \tilde{u}^t_m(s, X^{m,\vare}_s )ds
+\EE\int^t_0 \mathscr{L}^m_{1}(Y^{m,\vare}_s)\tilde{u}^t_m(s, X^{m,\vare}_s)ds\right|\nonumber\\
=\!\!\!\!\!\!\!\!&&\left|\EE\int^t_0 D_x \tilde{u}^t_m(s, X^{m,\vare}_s )\cdot [B^m(X^{m,\vare}_s,Y^{m,\vare}_s)-\bar{B}^m(X^{m,\vare}_s)] ds \right|
\end{eqnarray*}
Define $\rho(\vare):=\vare^{1-r}$ for $r\in (0,1)$. If $t\leq 2\rho(\vare)$, then by \eref{UE3} and the boundedness of $B$,
\begin{eqnarray}
\left|\EE\phi(X^{m,\vare}_{t})-\EE\phi(\bar{X}^m_{t})\right|\leq C_T\rho(\vare). \label{Small t}
\end{eqnarray}
If $t>2\rho(\vare)$,
\begin{eqnarray}
\left|\EE\phi(X^{m,\vare}_{t})-\EE\phi(\bar{X}^m_{t})\right|\leq\!\!\!\!\!\!\!\!&&\left|\EE\int^{t-\rho(\vare)}_{\rho(\vare)} D_x \tilde{u}^t_m(s, X^{m,\vare}_s )
\cdot [B^m(X^{m,\vare}_s,Y^{m,\vare}_s)-\bar{B}^m(X^{m,\vare}_s)]ds \right|\nonumber\\
&&+\left|\EE\int^{t}_{t-\rho(\vare)} D_x \tilde{u}^t_m(s, X^{m,\vare}_s )\cdot [B^m(X^{m,\vare}_s,Y^{m,\vare}_s)-\bar{B}^m(X^{m,\vare}_s)] ds \right|\nonumber\\
&&+\left|\EE\int^{\rho(\vare)}_0 D_x \tilde{u}^t_m(s, X^{m,\vare}_s )\cdot [B^m(X^{m,\vare}_s,Y^{m,\vare}_s)-\bar{B}^m(X^{m,\vare}_s)]ds \right|\nonumber\\
\leq\!\!\!\!\!\!\!\!&&\left|\EE\int^{t-\rho(\vare)}_{\rho(\vare)} D_x \tilde{u}^t_m(s, X^{m,\vare}_s )\cdot [B^m(X^{m,\vare}_s,Y^{m,\vare}_s)-\bar{B}^m(X^{m,\vare}_s)] ds \right|\nonumber\\
&&+C_T\rho(\vare).\label{F5.11}
\end{eqnarray}

For any $s\in [0,t], x,y\in H_m$, define
\begin{eqnarray*}
&&G^t_m(s,x,y):=D_x \tilde{u}^t_m(s, x)\cdot B^m(x,y)\\
&&\bar{G}^t_m(s,x):=\int_{H_m} G^t_m(s,x,y)\mu^{x,m}(dy)=D_x \tilde{u}^t_m(s, x)\cdot \bar{B}^m(x).
\end{eqnarray*}

By an argument similar to that used in the proof of Proposition \ref{P3.6}, we construct
$$
\tilde{\Phi}^t_m(s, x,y):=\int^{\infty}_0 \EE G^t_m(s,x,Y^{x,y,m}_r)-\bar{G}^t_m(s,x) dr, \quad s\in [0,t], x,y\in H_m,
$$
which is a solution of the following Poisson equation:
\begin{eqnarray}
-\mathscr{L}^m_{2}(x)\tilde{\Phi}^t_m(s,x,y)=G^t_m(s,x,y)-\bar{G}^t_m(s,x).\label{WPE}
\end{eqnarray}
Moreover, for any $T>0$, $t\in [0,T]$, $\delta\in (0,1]$, there exists $C_T,C_{T,\delta}>0$ such that the following estimates hold:
\begin{eqnarray}
&&\sup_{m\geq 1}|\partial_s \tilde{\Phi}^t_m(s,x,y)|\leq C_{T,\delta}(t-s)^{-1}(1+|x|)(1+|x|^{\delta}+|y|^{\delta}),\quad s\in (0,t];\label{E120}\\
&&\sup_{s\in [0, t],x\in H_m,m\geq 1}|\tilde{\Phi}^t_m(s,x,y)|\leq C_T(1+|x|+|y|);\label{E121}\\
&&\sup_{s\in [0, t],x\in H_m,m\geq 1}|D_x \tilde{\Phi}^t_m(s,x,y)\cdot h|\leq C_{T,\delta}(1+|x|^{\delta}+|y|^{\delta})|h|;\label{E122}\\
&&\sup_{s\in [0, t],x\in H_m,m\geq 1}|D_{xx}\tilde{\Phi}^t_m(s,x, y)\cdot(h,k)|\leq \!\!C_{T}(1+|x|+|y|)|h|\|k\|_{\kappa_1}. \label{E221}
\end{eqnarray}
We here only give the proof of \eref{E120}, and the proofs of \eref{E121}-\eref{E221} are omitted since it follows almost the same argument in Proposition \ref{P3.6}.

In fact, by \eref{ergodicity1}, \eref{UE2} and the boundedness of $B$, we have for small enough $\delta>0$,
\begin{eqnarray*}
|\partial_s \tilde{\Phi}^t_m(s,x,y)|\leq\!\!\!\!\!\!\!\!&&\int^{\infty}_{0} \partial_s D_x \tilde{u}^t_m(s, x)\cdot \left[\EE B^m(x,Y^{x,y,m}_r)-\bar{B}^m(x)\right] dr\\
\leq\!\!\!\!\!\!\!\!&&C_T(t-s)^{-1}(1+|x|)\int^{\infty}_{0}|\EE B^m(x,Y^{x,y,m}_r)-\bar{B}^m(x)|dr\\
\leq\!\!\!\!\!\!\!\!&&C_T(t-s)^{-1}(1+|x|)\int^{\infty}_{0}|\EE B^m(x,Y^{x,y,m}_r)-\bar{B}^m(x)|^{\delta}dr\\
\leq\!\!\!\!\!\!\!\!&& C_T(t-s)^{-1}(1+|x|)(1+|x|^{\delta}+|y|^{\delta})\int^{\infty}_{0}e^{-\frac{(\lambda_1-L_F)\delta r}{2}}dr\\
\leq\!\!\!\!\!\!\!\!&& C_{T,\delta}(t-s)^{-1}(1+|x|)(1+|x|^{\delta}+|y|^{\delta}).
\end{eqnarray*}


\vspace{0.2cm}
\textbf{Step 3.}Applying It\^o's formula and taking expectation, we get for any $t\in [2\rho(\vare),T]$,
\begin{eqnarray*}
&&\EE\tilde{\Phi}^t_m(t-\rho(\vare), X_{t-\rho(\vare)}^{m,\vare},Y^{m,\vare}_{t-\rho(\vare)})=\EE\tilde \Phi^t_m(\rho(\vare),X_{\rho(\vare)}^{m,\vare},Y^{m,\vare}_{\rho(\vare)})
+\EE\int^{t-\rho(\vare)}_{\rho(\vare)} \partial_s \tilde{\Phi}^t_m(s, X_{s}^{m,\vare},Y^{m,\vare}_{s})ds\\
&&+\EE\int^{t-\rho(\vare)}_{\rho(\vare)}\mathscr{L}^m_{1}(Y^{m,\vare}_{s})\tilde\Phi^t_m(s, X_{s}^{m,\vare},Y^{m,\vare}_{s})ds
+\frac{1}{\vare}\EE\int^{t-\rho(\vare)}_{\rho(\vare)} \mathscr{L}^m_{2}(X_{s}^{m,\vare})\tilde{\Phi}^t_m(s, X_{s}^{m,\vare},Y^{m,\vare}_{s})ds,
\end{eqnarray*}
which implies
\begin{eqnarray}
&&-\EE\int^{t-\rho(\vare)}_{\rho(\vare)} \mathscr{L}^m_{2}(X_{s}^{m,\vare})\tilde\Phi^t_m(s, X_{s}^{m,\vare},Y^{m,\vare}_{s})ds\nonumber\\
=\!\!\!\!\!\!\!\!&&\vare\big[\EE\tilde{\Phi}^t_m(\rho(\vare),X_{\rho(\vare)}^{m,\vare},Y^{m,\vare}_{\rho(\vare)})
-\EE\tilde{\Phi}^t_m(t-\rho(\vare), X_{t-\rho(\vare)}^{m,\vare},Y^{m,\vare}_{t-\rho(\vare)})
+\EE\int^{t-\rho(\vare)}_{\rho(\vare)} \partial_s \tilde{\Phi}^{t}_m(s, X_{s}^{m,\vare},Y^{m,\vare}_{s})ds\nonumber\\
&&+\EE\int^{t-\rho(\vare)}_{\rho(\vare)}\mathscr{L}^m_{1}(Y^{m,\vare}_{s})\tilde{\Phi}^{t}_m(s, X_{s}^{m,\vare},Y^{m,\vare}_{s})ds\big].\label{F3.39}
\end{eqnarray}

Combining  \eref{F5.11}, \eref{WPE} and \eref{F3.39}, we get for any $t\in [2\rho(\vare),T]$,
\begin{eqnarray*}
&&\left|\EE\phi(X^{m,\vare}_{t})-\EE\phi(\bar{X}^m_{t})\right|\leq C_T\rho(\vare)+\left|\EE\int^{t-\rho(\vare)}_{\rho(\vare)} \mathscr{L}^m_{2}(X_{s}^{m,\vare})\tilde{\Phi}^t_m(s, X_{s}^{m,\vare},Y^{m,\vare}_{s})ds\right|\\
\leq\!\!\!\!\!\!\!\!&&C_T\rho(\vare)+\vare\Bigg[\EE|\tilde{\Phi}^t_m(\rho(\vare),X_{\rho(\vare)}^{m,\vare},Y^{m,\vare}_{\rho(\vare)})|+\left|\EE\tilde{\Phi}^t_m(t-\rho(\vare), X_{t-\rho(\vare)}^{m,\vare},Y^{m,\vare}_{t-\rho(\vare)})\right|\nonumber\\
&&+\EE\int^{t-\rho(\vare)}_{\rho(\vare)} \left|\partial_s \tilde{\Phi}^t_m(s, X_{s}^{m,\vare},Y^{m,\vare}_{s})\right|ds+\EE\int^{t-\rho(\vare)}_{\rho(\vare)}\left|\mathscr{L}^m_{1}(Y^{m,\vare}_{s})\tilde{\Phi}^t_m(s, X_{s}^{m,\vare},Y^{m,\vare}_{s})\right|ds\Bigg]\nonumber\\
:=\!\!\!\!\!\!\!\!&&C_T\rho(\vare)+\vare\sum^4_{k=1}\tilde{\Lambda}^{m,\vare}_k(t).
\end{eqnarray*}

For the terms $\tilde{\Lambda}^{m,\vare}_1(t)$ and  $\tilde{\Lambda}^{m,\vare}_2(t)$. By estimates \eref{E121} and \eref{Yvare}, it is easy to get
\begin{eqnarray}
\tilde{\Lambda}^{m,\vare}_1(t)+\tilde{\Lambda}^{m,\vare}_2(t)\leq C_{T}\sup_{t\in [0,T]}\EE(1+|X^{m,\vare}_{t}|+|Y^{m,\vare}_{t}|)\leq C_T(1+|x|+|y|).\label{Lambda12(T)}
\end{eqnarray}

For the terms $\tilde{\Lambda}^{m,\vare}_3(t)$. By \eref{E120}, \eref{Xvare} and \eref{Yvare}, we have
\begin{eqnarray}
\tilde{\Lambda}^{m,\vare}_3(t)\leq\!\!\!\!\!\!\!\!&& C_{T,\delta}\int^{t-\rho(\vare)}_0(t-s)^{-1}\EE\left[(1+|X^{m,\vare}_{s}|)(1+|X^{m,\vare}_{s}|^{\delta}+|Y^{m,\vare}_{s}|^{\delta})\right]ds\nonumber\\
\leq\!\!\!\!\!\!\!\!&&C_{T,\delta} \ln\left(\frac{T}{\rho(\vare)}\right)(1+|x|^{1+\delta}+|y|^{1+\delta}).\label{Lambda3(T)}
\end{eqnarray}

For the terms $\tilde{\Lambda}^{m,\vare}_4(t)$. It is easy to see
\begin{eqnarray}
\tilde{\Lambda}^{m,\vare}_4(t)\leq\!\!\!\!\!\!\!\!&&\int^t_{\rho(\vare)} \EE\left|D_x\tilde{\Phi}^t_m(X_{s}^{m,\vare},Y^{m,\vare}_{s})\cdot AX_{s}^{m,\vare} \right|ds\nonumber\\
&&+\int^t_{\rho(\vare)} \EE\left|D_x\tilde{\Phi}^t_m(X_{s}^{m,\vare},Y^{m,\vare}_{s})\cdot B^m(X_{s}^{m,\vare},Y^{m,\vare}_{s}) \right|ds\nonumber\\
&&+\int^t_0\EE\left|\sum^m_{k=1}\beta^{\alpha}_k\int_{\RR}\big[\tilde{\Phi}^t_m(X_{s}^{m,\vare}+e_k z,Y^{m,\vare}_{s})-\tilde{\Phi}^t_m(X_{s}^{m,\vare},Y^{m,\vare}_{s})\right.\nonumber\\
&&\quad\quad\left.-D_x\tilde{\Phi}^t_m(X_{s}^{m,\vare},Y^{m,\vare}_{s})\cdot( e_k z) 1_{\{|z|\leq 1\}}\big]\nu (dz)\right|ds
:=\sum^3_{i=1}\tilde{\Lambda}^{m,\vare}_{4i}(t).\label{Lambda4(T)}
\end{eqnarray}

By \eref{E122} and \eref{Xvare2}, it is easy to see for any $r\in (0,1)$
\begin{eqnarray}
\tilde{\Lambda}^{m,\vare}_{41}(t)\leq\!\!\!\!\!\!\!\!&&C_{T,\delta}\EE\int^T_{\rho(\vare)} \|X_{s}^{m,\vare}\|_{2}(1+|X_{s}^{m,\vare}|^{\delta}+|Y^{m,\vare}_{s}|^{\delta})ds\nonumber\\
\leq\!\!\!\!\!\!\!\!&&C_{T,\delta}\int^T_{\rho(\vare)} \left[\EE\|X_{s}^{m,\vare}\|^p_{2}\right]^{1/p}\left[\EE(1+|X_{s}^{m,\vare}|^{\frac{\delta p}{p-1}}+|Y^{m,\vare}_{s}|^{\frac{\delta p}{p-1}})\right]^{\frac{p-1}{p}}ds\nonumber\\
\leq\!\!\!\!\!\!\!\!&&C_{T}\int^{T}_{\rho(\vare)} (s^{-1}+\vare^{-r})(1+|x|^{1+\delta}+|y|^{1+\delta})ds\nonumber\\
\leq\!\!\!\!\!\!\!\!&&C_T \left[\ln\left(\frac{T}{\rho(\vare)}\right)+\vare^{-r}\right](1+|x|^{1+\delta}+|y|^{1+\delta}).
\end{eqnarray}

By \eref{E122}, \eref{E221}, \eref{Xvare} and \eref{Yvare}, we have
\begin{eqnarray}
\tilde{\Lambda}^{m,\vare}_{42}(t)\leq\!\!\!\!\!\!\!\!&&C_{T}\EE\int^T_0 (1+|X_{s}^{m,\vare}|+|Y^{m,\vare}_{s}|)(1+|Y^{m,\vare}_{s}|^{\delta})ds\nonumber\\
\leq\!\!\!\!\!\!\!\!&&C_{T}(1+|x|^{1+\delta}+|y|^{1+\delta}).
\end{eqnarray}

By following the same argument in the estimating $\Lambda^{m,\vare}_{53}(T)$, we have for any $\delta<\alpha-1$,
\begin{eqnarray}
\tilde{\Lambda}^{m,\vare}_{43}(t)\leq C_{\delta,T}(1+|x|+|y|).\label{F5.22}
\end{eqnarray}

Finally, combining estimates \eref{Small t}, \eref{F5.11}, \eref{Lambda12(T)}-\eref{F5.22},  we final obtain
\begin{eqnarray*}
&&\sup_{t\in[0,T],m\geq 1}\left|\EE\phi(X^{m,\vare}_{t})-\EE\phi(\bar{X}^m_{t})\right|\\
\leq\!\!\!\!\!\!\!\!&&C_T\vare^{1-r}(1+|x|^{1+\delta}+|y|^{1+\delta})+C_T\vare \ln\left(\frac{T}{\vare^{1-r}}\right)(1+|x|^{1+\delta}+|y|^{1+\delta})\\
\leq\!\!\!\!\!\!\!\!&&C_{T,r}\vare^{1-r}(1+|x|^{1+\delta}+|y|^{1+\delta}).
\end{eqnarray*}
The proof is complete.

\section{Appendix}

In this section, we give some a priori estimates of the solution $(X_{t}^{\varepsilon}, Y_{t}^{\varepsilon})$ (see Lemma \ref{L6.1}), which is used to study the Galerkin approximation of the system \eref{main equation} (see Lemma \ref{GA1}). Then we study the increment of the time of solution $(X_{t}^{\varepsilon}, Y_{t}^{\varepsilon})$ (see Lemma \ref{L6.3}). Finally, the finite dimensional approximation of the frozen equation \eref{FZE} is given (see Lemma \ref{GA2}).

\begin{lemma} \label{L6.1}
For any $x,y\in H$, $1\leq p<\alpha$ and $T>0$, there exist constants $C_{p,T},\ C_{T}>0$ such that the solution $(X^{\vare}_t, Y^{\vare}_t)$ of system \eref{main equation} satisfies
\begin{eqnarray}
\sup_{t\in [0,T]}\mathbb{E}|X_{t}^{\varepsilon}|^p \leq  C_{p,T}(1+ |x|^p+|y|^p),\quad \forall \vare>0,\label{AXvare}
\end{eqnarray}
\begin{eqnarray}
\sup_{t\in [0,T]}\mathbb{E} |Y_{t}^{\varepsilon} |^p\leq C_{p}(1+|x|^p+|y|^p ),\quad \forall \vare>0.\label{AYvare}
\end{eqnarray}
\end{lemma}

\begin{proof}

Define $\tilde Z_t:=\frac{1}{\vare^{1/ \alpha}}Z_{t\vare}$, which is also a cylindrical $\alpha$-stable process. Then by \cite[(4.12)]{PZ}, for any $p\in (1,\alpha)$,
\begin{align*}
\EE\left|\frac{1}{\vare^{1/\alpha}}\int^t_0 e^{(t-s)A/\varepsilon}dZ_s\right|^p=&\EE\left|\int^{t/\vare}_0 e^{(t/\vare-s)A}d\tilde Z_s\right|^{p}\\
\leq&C\left(\sum^{\infty}_{k=1}\gamma^{\alpha}_k \frac{1-e^{-\alpha \lambda_k t/\vare}}{\alpha \lambda_k} \right)^{p/\alpha}\\
\leq&C\left(\sum^{\infty}_{k=1}\frac{\gamma^{\alpha}_k }{\alpha \lambda_k} \right)^{p/\alpha}<\infty.
\end{align*}
Then by Minkowski's inequality, we get for any $p\in (1,\alpha)$ and $0<t\leq T$,
\begin{eqnarray*}
\left(\EE|Y^{\varepsilon}_t|^p\right)^{1/p}
\leq\!\!\!\!\!\!\!\!&&|e^{tA/\vare}y|+\left[\EE\left(\frac{1}{\vare}\int^t_0|e^{\frac{(t-s)A}{\vare}}F(X^{\vare}_s,Y^{\vare}_s)|ds\right)^p\right]^{1/p}\nonumber\\
&&+\left[\EE\left|\frac{1}{\vare^{1/\alpha}}\int^t_0 e^{(t-s)A/\varepsilon}dZ_s\right|^p\right]^{1/p}\\
\leq\!\!\!\!\!\!\!\!&&|y|+\frac{1}{\vare}\int^t_0|e^{\frac{-\lambda_1(t-s)}{\vare}}\left[C+C\left(\EE|X^{\vare}_s|^p\right)^{1/p}+L_F\left(\EE|Y^{\vare}_s|^p\right)^{1/p}\right]\nonumber\\
&&+C\left(\sum^{\infty}_{k=1}\frac{\gamma^{\alpha}_k }{\alpha \lambda_k} \right)^{1/\alpha}\\
\leq\!\!\!\!\!\!\!\!&&C(1+|y|)+\frac{L_{F}}{\lambda_1}\sup_{0\leq t\leq T}\left(\EE|Y^{\varepsilon}_t|^p\right)^{1/p}+C\sup_{t\in [0,T]}\left(\mathbb{E}|X_{t}^{\varepsilon} |^p\right)^{1/p}.
\end{eqnarray*}
Thus by $L_F<\lambda_1$, it follows
\begin{eqnarray}
\sup_{t\in [0,T]}\EE|Y^{m,\varepsilon}_t|^p
\leq C_p(1+|y|^p)+C_p\sup_{t\in [0,T]}\mathbb{E}|X_{t}^{\varepsilon} |^p.\label{F6.3}
\end{eqnarray}

For any $1<p<\alpha$, by \eref{F6.3} and \eref{LA} in Remark \ref{Re2}, we obtain that
\begin{eqnarray*}
\sup_{t\in [0,T]}\mathbb{E}|X_{t}^{\varepsilon} |^p\leq\!\!\!\!\!\!\!\!&&C_p|x|^p\!+\!C_p\int^T_0 \!\EE|X^{\varepsilon}_s|^p ds\!+\!C_p\int^T_0\!\EE|Y^{\varepsilon}_s|^pds\!+\!C_p\sup_{t\in [0,T]}\mathbb{E}\left|\int^t_0 e^{(t-s)A}dL_s\right|^p\\
\leq\!\!\!\!\!\!\!\!&&C_{p}(1+|x|^p+|y|^p)+C_p\int^T_0 \sup_{r\leq s}\EE|X^{\varepsilon}_r|^p ds.
\end{eqnarray*}
Then by Gronwall's inequality, we have
\begin{eqnarray*}
\sup_{t\in [0,T]}\mathbb{E}|X_{t}^{\varepsilon}|^p\leq\!\!\!\!\!\!\!\!&&C_{p,T}(1+|x|^p+|y|^p),
\end{eqnarray*}
which also implies that \eref{AYvare} holds easily. The proof is complete.
\end{proof}

\vspace{0.3cm}
Recall the Galerkin approximation \eref{Ga mainE} of system \eref{main equation}. We have the following approximation.

\begin{lemma} \label{GA1}
For any $\vare>0, (x,y)\in H\times H$ and $m\in \mathbb{N}_{+}$, system \eref{Ga mainE} has a unique mild solution $(X^{m,\vare}_t, Y^{m,\vare}_t)\in H\times H$, i.e., $\PP$-a.s.,
\begin{equation}\left\{\begin{array}{l}\label{A mild solution of F}
\displaystyle
X^{m,\varepsilon}_t=e^{tA}x^m+\int^t_0e^{(t-s)A}B^m(X^{m,\varepsilon}_s, Y^{m,\varepsilon}_s)ds+\int^t_0 e^{(t-s)A}d\bar {L}^m_s,\nonumber\vspace{2mm}\\
Y^{m,\varepsilon}_t=e^{tA/\varepsilon}y^m+\frac{1}{\varepsilon}\int^t_0e^{(t-s)A/\varepsilon}F^m(X^{m,\varepsilon}_s,Y^{m,\varepsilon}_s)ds
+\frac{1}{\vare^{1/\alpha}}\int^t_0 e^{(t-s)A/\varepsilon}d\bar {Z}^m_s.\nonumber
\end{array}\right.
\end{equation}
Moreover, for any $1\leq p<\alpha$ and $T>0$, there exist constants $C_{p,T},\ C_{T}>0$ such that for any $\vare>0$,
\begin{eqnarray}
&&\sup_{m\geq 1, t\in [0,T]}\mathbb{E}|X_{t}^{m,\varepsilon}|^p \leq  C_{p,T}(1+ |x|^p+|y|^p);\label{Xvare}\\
&&\sup_{m\geq 1, t\in [0,T]}\mathbb{E}|Y_{t}^{m,\varepsilon} |^p\leq C_{p,T}(1+|x|^p+|y|^p);\label{Yvare}\\
&&\lim_{m\rightarrow \infty}\sup_{t\in[0,T]}\EE|X^{m,\vare}_t-X^{\vare}_t|^p=0, \quad \forall\vare>0.\label{FA1}
\end{eqnarray}
\end{lemma}

\begin{proof}
Under the assumptions \ref{A1} and \ref{A2}, it is easy to show that the existence and uniqueness of the mild solution of system \eref{Ga mainE} . The estimates \eref{Xvare} and \eref{Yvare} can be proved by following the same argument as in the proof of Lemma \ref{L6.1}. Next, we prove the approximation \eref{FA1}.
It is easy to see that for any $t>0$,
\begin{eqnarray*}
X^{m,\vare}_t-X^{\vare}_t=\!\!\!\!\!\!\!\!&&e^{tA}(x^m-x)+\int^{t}_{0}e^{(t-s)A}(\pi_m-I)B(X^{\vare}_s, Y^{\vare}_s)ds\nonumber\\
&&\!\!\!\!\!\!\!\!+\int^{t}_{0}e^{(t-s)A}\left[B^{m}(X^{m,\vare}_s,Y^{m,\vare}_s)-B^{m}(X^{\vare}_s,Y^{\vare}_s)\right]ds\!+\!\left[\int^t_0 e^{(t-s)A}d \bar L^{m}_s-\int^t_0 e^{(t-s)A}dL_s\right].
\end{eqnarray*}
Then for any $T>0$, $p\in [1,\alpha)$, we get
\begin{eqnarray}
\sup_{t\in [0,T]}\EE|X^{m,\vare}_t-X^{\vare}_t|^p\leq\!\!\!\!\!\!\!\!&&C_p|x^m-x|^p+C_p\int^{T}_{0}\EE|(\pi_m-I)B(X^{\vare}_s, Y^{\vare}_s)|^pds\nonumber\\
&&+C_p\EE\left[\int^{T}_{0}|B^{m}(X^{m,\vare}_s,Y^{m,\vare}_s)-B^{m}(X^{\vare}_s,Y^{\vare}_s)|ds\right]^p\nonumber\\
&&+C_p\sup_{t\in [0,T]}\EE\left|\int^t_0 e^{(t-s)A}d \bar L^{m}_s-\int^t_0 e^{(t-s)A}dL_s\right|^p\nonumber\\
\leq\!\!\!\!\!\!\!\!&&C_p|x^m-x|^p+C_p\int^{T}_{0}\EE|(\pi_m-I)B(X^{\vare}_s, Y^{\vare}_s)|^pds\nonumber\\
&&+C_p\EE\left(\int^{T}_{0}|X^{m,\vare}_s-X^{\vare}_s|+|Y^{m,\vare}_s-Y^{\vare}_s|ds\right)^p\nonumber\\
&&+C_p\sup_{t\in [0,T]}\EE\left|\int^t_0 e^{(t-s)A}d \bar L^{m}_s-\int^t_0 e^{(t-s)A}dL_s\right|^p.\label{X^m-X}
\end{eqnarray}
On the other hand,
\begin{eqnarray*}
Y^{m,\vare}_t-Y^{\vare}_t=\!\!\!\!\!\!\!\!&&e^{tA/\vare}(y^m-y)+\frac{1}{\vare}\int^{t}_{0}e^{(t-s)A/\vare}(\pi_m-I)F(X^{\vare}_s, Y^{\vare}_s)ds \nonumber\\ &&+\frac{1}{\vare}\int^{t}_{0}e^{(t-s)A/\vare}[F^{m}(X^{m,\vare}_s,Y^{m,\vare}_s)-F^{m}(X^{\vare}_s,Y^{\vare}_s)]ds\\
&&+\frac{1}{\vare^{1/\alpha}}\int^t_0 e^{(t-s)A/\vare}d \bar{Z}^m_s-\frac{1}{\vare^{1/\alpha}}\int^t_0 e^{(t-s)A/\vare}d Z_s.
\end{eqnarray*}
It is clear that for any $T>0$,
\begin{eqnarray*}
\int^T_0|Y^{m,\vare}_t-Y^{\vare}_t|dt\leq\!\!\!\!\!\!\!\!&&C|y^m-y|+\frac{1}{\vare}\int^{T}_{0}\int^{t}_{0}e^{-\frac{\lambda_1(t-s)}{\vare}}|(\pi_m-I)F(X^{\vare}_s, Y^{\vare}_s)|dsdt\nonumber\\
&&+\frac{1}{\vare}\int^T_0\int^{t}_{0}e^{-\frac{\lambda_1(t-s)}{\vare}}|F^{m}(X^{m,\vare}_s,Y^{m,\vare}_s)-F^{m}(X^{\vare}_s,Y^{\vare}_s)|dsdt\\
&&+\int^T_0\left|\frac{1}{\vare^{1/\alpha}}\int^t_0 e^{(t-s)A/\vare}d \bar{Z}^m_s-\frac{1}{\vare^{1/\alpha}}\int^t_0 e^{(t-s)A/\vare}d Z_s\right|dt\\
\leq\!\!\!\!\!\!\!\!&&C_T|y^m-y|+\frac{1}{\vare}\int^{T}_{0}|(\pi_m-I)F(X^{\vare}_s, Y^{\vare}_s)|\int^{T}_{s}e^{-\frac{\lambda_1(t-s)}{\vare}}dtds\nonumber\\
&&+\frac{1}{\vare}\int^T_0\left(C|X^{m,\vare}_s-X^{\vare}_s|+L_F|Y^{m,\vare}_s-Y^{\vare}_s|\right)\int^{T}_{s}e^{-\frac{\lambda_1(t-s)}{\vare}}dtds\\
&&+\int^T_0\left|\frac{1}{\vare^{1/\alpha}}\int^t_0 e^{(t-s)A/\vare}d \bar{Z}^m_s-\frac{1}{\vare^{1/\alpha}}\int^t_0 e^{(t-s)A/\vare}d Z_s\right|dt\\
\leq\!\!\!\!\!\!\!\!&&C|y^m-y|+\frac{1}{\lambda_1}\int^{T}_{0}|(\pi_m-I)F(X^{\vare}_s, Y^{\vare}_s)|ds\nonumber\\
&&+\frac{L_F}{\lambda_1}\int^T_0|Y^{m,\vare}_s-Y^{\vare}_s|ds+\frac{C}{\lambda_1}\int^T_0|X^{m,\vare}_s-X^{\vare}_s|ds\\
&&+\int^T_0\left|\frac{1}{\vare^{1/\alpha}}\int^t_0 e^{(t-s)A/\vare}d \bar{Z}^m_s-\frac{1}{\vare^{1/\alpha}}\int^t_0 e^{(t-s)A/\vare}d Z_s\right|dt.
\end{eqnarray*}
By the condition $\lambda_1>L_F$ in assumption \ref{A2}, it follows
\begin{eqnarray}
\int^T_0|Y^{m,\vare}_t-Y^{\vare}_t|dt
\leq\!\!\!\!\!\!\!\!&&C|y^m-y|+\frac{C}{\lambda_1}\int^{T}_{0}|(\pi_m-I)F(X^{\vare}_s, Y^{\vare}_s)|ds+\frac{C}{\lambda_1}\int^T_0|X^{m,\vare}_s-X^{\vare}_s|ds\nonumber\\
&&+C\int^T_0\left|\frac{1}{\vare^{1/\alpha}}\int^t_0 e^{(t-s)A/\vare}d \bar{Z}^m_s-\frac{1}{\vare^{1/\alpha}}\int^t_0 e^{(t-s)A/\vare}d Z_s\right|dt.\label{Y^m-Y}
\end{eqnarray}
Then by \eref{X^m-X} and \eref{Y^m-Y}, we obtain
\begin{eqnarray*}
\sup_{t\in [0,T]}\EE|X^{m,\vare}_t-X^{\vare}_t|^p\leq\!\!\!\!\!\!\!\!&&C_p|x^m-x|^p+C_{p}|y^m-y|^p+C_{p,T}\EE\int^{T}_{0}|X^{m,\vare}_s-X^{\vare}_s|^p ds\nonumber\\
&&\!\!\!\!\!\!\!\!+C_{p,T}\!\int^{T}_{0}\!\!\EE|(\pi_m-I)B(X^{\vare}_s, Y^{\vare}_s)|^pds\!+\!C_{p,T}\!\int^{T}_{0}\!\!\EE|(\pi_m-I)F(X^{\vare}_s, Y^{\vare}_s)|^p ds\nonumber\\
&&\!\!\!\!\!\!\!\!+C_{p,T}\int^T_0\EE\left|\frac{1}{\vare^{1/\alpha}}\int^t_0 e^{(t-s)A/\vare}d \bar{Z}^m_s-\frac{1}{\vare^{1/\alpha}}\int^t_0 e^{(t-s)A/\vare}d Z_s\right|^pdt\nonumber\\
&&\!\!\!\!\!\!\!\!+C_p\sup_{t\in [0,T]}\EE\left|\int^t_0 e^{(t-s)A}d \bar L^{m}_s-\int^t_0 e^{(t-s)A}dL_s\right|^p.
\end{eqnarray*}
The Gronwall's inequality implies
\begin{eqnarray*}
&&\sup_{t\in [0,T]}\EE|X^{m,\vare}_t-X^{\vare}_t|^p\\
\leq\!\!\!\!\!\!\!\!&&C_{p,T}(|x^m-x|^p+|y^m-y|^p)+C_{p,T}\int^{T}_{0}\EE|(\pi_m-I)B(X^{\vare}_s, Y^{\vare}_s)|^pds\nonumber\\
&&\!\!\!\!\!\!\!\!+C_{p,T}\int^{T}_{0}\!\!\EE|(\pi_m-I)F(X^{\vare}_s, Y^{\vare}_s)|^p ds\!+\!C_{p,T}\!\int^T_0\EE\left|\frac{1}{\vare^{1/\alpha}}\int^t_0 \!\!e^{(t-s)A/\vare}d \bar{Z}^m_s\!-\!\frac{1}{\vare^{1/\alpha}}\int^t_0\!\! e^{(t-s)A/\vare}d Z_s\right|^pdt\nonumber\\
&&\!\!\!\!\!\!\!\!+C_{p,T}\sup_{t\in [0,T]}\EE\left|\int^t_0 e^{(t-s)A}d \bar L^{m}_s-\int^t_0 e^{(t-s)A}dL_s\right|^p.
\end{eqnarray*}

It is clear that as $m\rightarrow\infty$,
$$
|x^m-x|^p\rightarrow 0,\quad |y^m-y|^p\rightarrow 0.
$$
By the a prior estimate of $(X^{\vare}_s, Y^{\vare}_s)$ and the dominated convergence theorem,
\begin{eqnarray*}
&&\lim_{m\rightarrow \infty}\int^{T}_{0}\EE|(\pi_m-I)B(X^{\vare}_s, Y^{\vare}_s)|^pds=\int^{T}_{0}\EE\lim_{m\rightarrow \infty}|(\pi_m-I)B(X^{\vare}_s, Y^{\vare}_s)|^pds=0,\\
&&\lim_{m\rightarrow \infty}\int^{T}_{0}\EE|(\pi_m-I)F(X^{\vare}_s, Y^{\vare}_s)|^p ds=\int^{T}_{0}\EE\lim_{m\rightarrow \infty}|(\pi_m-I)F(X^{\vare}_s, Y^{\vare}_s)|^p ds=0.
\end{eqnarray*}
By assumption \ref{A2} and Remark \ref{Re2}, as $m\rightarrow\infty$,
\begin{eqnarray*}
&&\sup_{t\in [0,T]}\EE\left|\int^t_0 e^{(t-s)A}d \bar L^{m}_s-\int^t_0 e^{(t-s)A}dL_s\right|^p\leq C_{p}\left(\sum^{\infty}_{k=m+1}\frac{\beta^{\alpha}_k}{\lambda_k}\right)^{p/\alpha}\!\!\!\rightarrow 0,\\
&&\int^T_0\EE\left|\frac{1}{\vare^{1/\alpha}}\int^t_0 e^{(t-s)A/\vare}d \bar{Z}^m_s-\frac{1}{\vare^{1/\alpha}}\int^t_0 e^{(t-s)A/\vare}d Z_s\right|^pdt\leq C_{p,T}\left(\sum^{\infty}_{k=m+1}\frac{\gamma^{\alpha}_k}{\lambda_k}\right)^{p/\alpha}\!\!\!\rightarrow 0.
\end{eqnarray*}
Hence by the discussion above,  we final get
\begin{eqnarray*}
\lim_{m\rightarrow \infty}\sup_{t\in [0,T]}\EE|X^{m,\vare}_t-X^{\vare}_t|^p=0.
\end{eqnarray*}
The proof is complete.
\end{proof}

\begin{lemma} \label{L6.3}
For any $ (x,y)\in H\times H$, $m\in \mathbb{N}_{+}$, $p\in [1,\alpha)$ and $T>0$, there exist constants $C_{p,T}, C_{T}>0$ such that any  $\vare\in (0,1], \eta\in (0,1)$ and $0<s\leq t<T$,
\begin{eqnarray}
&&\sup_{m\geq 1}\left(\EE|X^{m,\vare}_t-X^{m,\vare}_s|^p\right)^{1/p}\leq C_T(t-s)^{\frac{\eta}{2}}s^{-\frac{\eta}{2}}(1+|x|+|y|),\quad\label{THXvare}\\
&&\sup_{m\geq 1}\left(\EE|Y^{m,\vare}_t-Y^{m,\vare}_s|^p\right)^{1/p}\leq C_T\left(\frac{t-s}{\vare}\right)^{\frac{\eta}{2}}s^{-\frac{\eta}{2}}(1+|x|+|y|).\label{THYvare}
\end{eqnarray}
\end{lemma}
\begin{proof}
By  \eref{P3}, \eref{LA} and Minkowski's inequality, we get for any $p\in (1,\alpha)$, $\eta\in (0,1)$ and $0<t\leq T$,
\begin{eqnarray}
\EE\|X^{m,\varepsilon}_t\|^p_{\eta}
\leq\!\!\!\!\!\!\!\!&&C_p\|e^{tA}x\|^p_{\eta}+C_p\EE\left(\int^t_0\|e^{(t-s)A}B^m(X^{m,\vare}_s,Y^{m,\vare}_s)\|_{\eta}ds\right)^p
+C_p\EE\left\|\int^t_0 e^{(t-s)A}d\bar{L}^m_s\right\|^p_{\eta}\nonumber\\
\leq\!\!\!\!\!\!\!\!&&C_p t^{-\frac{p\eta}{2}}|x|^p+C_p\left[\int^t_0(t-s)^{-\eta/2}\left[\EE(1+|X^{m,\varepsilon}_s|^p+|Y^{m,\varepsilon}_s|^p)\right]^{1/p}ds\right]^p+C_p\nonumber\\
\leq\!\!\!\!\!\!\!\!&&C_{p,T} t^{-\frac{p\eta}{2}}(1+|x|^p+|y|^p).\label{X_theta}
\end{eqnarray}

Note that
$$
X^{m,\varepsilon}_t=e^{(t-s)A}X^{m,\varepsilon}_s+\int^t_s e^{(t-r)A}B^m(X^{m,\varepsilon}_r, Y^{m,\varepsilon}_r)dr+\int^t_s e^{(t-r)A}d\bar {L}^{m}_r
$$
It follows from \eref{P4} that for any $0<s\leq t\leq T$,
\begin{eqnarray*}
\EE|X^{m,\varepsilon}_t-X^{m,\varepsilon}_s|^p\leq\!\!\!\!\!\!\!\!&& C_p(t-s)^{\frac{p\eta}{2}}\EE\|X^{m,\varepsilon}_s\|^p_{\eta}+C_p\left(\int^t_s \left[\EE|B^m(X^{m,\varepsilon}_r, Y^{m,\varepsilon}_r)|^p\right]^{1/p}dr\right)^p\\
&&+C_p\left[\sum^{\infty}_{k=1}\frac{\beta^{\alpha}_k(1-e^{-\lambda_k (t-s)})}{\lambda_k}\right]^{p/\alpha}\\
\leq\!\!\!\!\!\!\!\!&&C_{p,T}(t-s)^{\frac{p\eta}{2}}s^{-\frac{p\eta}{2}}(1+|x|^p+|y|^p)+C(t-s)^p(1+|x|^p+|y|^p)\\
&&+C_p(t-s)^{\frac{p\eta}{2}}\left[\sum^{\infty}_{k=1}\frac{\beta^{\alpha}_k}{\lambda^{1-\eta\alpha/2}_k}\right]^{p/\alpha}\\
\leq\!\!\!\!\!\!\!\!&&C_{p,T}(t-s)^{\frac{p\eta}{2}}s^{-\frac{p\eta}{2}}(1+|x|^p+|y|^p).
\end{eqnarray*}

Define $\tilde Z_t:=\frac{1}{\vare^{1/ \alpha}}Z_{t\vare}$, which is also a cylindrical $\alpha$-stable process. Then by \eref{LA} , for any $\eta\in(0,1)$ we have
\begin{align*}
\EE\left\|\frac{1}{\vare^{1/\alpha}}\int^t_0 e^{(t-s)A/\varepsilon}dZ_s\right\|^p_{\eta}=&\EE\left\|\int^{t/\vare}_0 e^{(t/\vare-s)A}d\tilde Z_s\right\|^p_{\eta}\\
\leq&C\left(\sum_{k}\gamma^{\alpha}_k \frac{1-e^{-\alpha \lambda_k t/\vare}}{\alpha \lambda^{1-\alpha\eta/2}_k} \right)^{p/\alpha}\\
\leq&C\left(\sum_{k}\frac{\gamma^{\alpha}_k }{\alpha \lambda^{1-\alpha\eta/2}_k} \right)^{p/\alpha}\\
\leq&C_{\alpha}\left(\sum_{k}\gamma^{\alpha}_k \right)^{p/\alpha}\leq C_{\alpha,p}.
\end{align*}
Similarly, by  \eref{P3} and Minkowski's inequality, we get for any $\eta\in (0,1)$ and $0<t\leq T$,
\begin{eqnarray}
\EE\|Y^{m,\varepsilon}_t\|^p_{\eta}
\leq\!\!\!\!\!\!\!\!&&C_p\|e^{tA/\vare}y\|^p_{\eta}+C_p\EE\left(\frac{1}{\vare}\int^t_0\|e^{\frac{(t-s)A}{\vare}}F^m(X^{m,\vare}_s,Y^{m,\vare}_s)\|_{\eta}ds\right)^p\nonumber\\
&&+C_p\EE\left\|\frac{1}{\vare^{1/\alpha}}\int^t_0 e^{\frac{(t-s)A}{\vare}}d\bar{Z}^{m}_s\right\|^p_{\eta}\nonumber\\
\leq\!\!\!\!\!\!\!\!&&C_p\left(\frac{t}{\vare}\right)^{-\frac{p\eta}{2}}|y|^p+C\left[\frac{1}{\vare}\int^t_0\left(\frac{t-s}{\vare}\right)^{-\frac{\eta}{2}}e^{-\frac{(t-s)\lambda_1}{2\vare}}\left[\EE(1+|X^{m,\varepsilon}_s|^p+|Y^{m,\varepsilon}_s|^p)\right]^{1/p}ds\right]^p\nonumber\\
&&+C_{\alpha,p}\nonumber\\
\leq\!\!\!\!\!\!\!\!&&C_T t^{-\frac{p\eta}{2}}(1+|x|^p+|y|^p).\label{Y_theta}
\end{eqnarray}

Note that
\begin{eqnarray*}
Y^{m,\varepsilon}_t=e^{\frac{(t-s)A}{\vare}}Y^{m,\varepsilon}_s+\frac{1}{\vare}\int^t_s e^{\frac{(t-r)A}{\vare}}F^m(X^{m,\varepsilon}_r, Y^{m,\varepsilon}_r)dr+\frac{1}{\vare^{1/\alpha}}\int^t_s e^{\frac{(t-r)A}{\vare}}d\bar {Z}^{m}_r.
\end{eqnarray*}

Thus by \eref{P4}, it follows that for any $0<s<t\leq T$,
\begin{eqnarray*}
&&\EE|Y^{m,\varepsilon}_t-Y^{m,\varepsilon}_s|^p\\
\leq\!\!\!\!\!\!\!\!&&C_p\EE|e^{\frac{(t-s)A}{\vare}}Y^{m,\varepsilon}_s-Y^{m,\varepsilon}_s|^p +C_p\EE\left[\frac{1}{\vare}\int^t_s \left|e^{\frac{(t-r)A}{\vare}}F^m_2(X^{m,\varepsilon}_r, Y^{m,\varepsilon}_r)\right|dr\right]^p\\
&&+C_p\EE\left|\frac{1}{\vare^{1/\alpha}}\int^t_s e^{\frac{(t-r)A}{\vare}}d\bar {Z}^{m}_r\right|^p\\
\leq\!\!\!\!\!\!\!\!&&C_p\left(\frac{t-s}{\vare}\right)^{\frac{p\eta}{2}}\EE\|Y^{m,\varepsilon}_s\|^p_{\eta}+C_p\left[\frac{1}{\vare}\int^t_s e^{-\frac{(t-r)\lambda_1}{\vare}}
\left[\EE(1+|X^{m,\varepsilon}_r|^p+|Y^{m,\varepsilon}_r|^p)\right]^{1/p}dr\right]^p\nonumber\\
&&+C_p\left[\sum^{\infty}_{k=1} \frac{\gamma^{\alpha}_k(1-e^{-\alpha \lambda_k(t-s)/\vare})}{\lambda_k}\right]^{p/\alpha}\nonumber\\
\leq\!\!\!\!\!\!\!\!&&C_{p,T}\left(\frac{t-s}{\vare}\right)^{\frac{p\eta}{2}}s^{-\frac{p\eta}{2}}(1+|x|^p+|y|^p)+C_{p,T}\left(\frac{t-s}{\vare}\right)^{\frac{p\eta}{2}}(1+|x|^p+|y|^p)\\
&&+C_{p,T}\left(\frac{t-s}{\vare}\right)^{\frac{p\eta}{2}}\left(\sum^{\infty}_{k=1} \frac{\gamma^{\alpha}_k}{\lambda^{1-\alpha\eta/2}_k}\right)^{p/\alpha}\\
\leq\!\!\!\!\!\!\!\!&&C_{p,T}\left(\frac{t-s}{\vare}\right)^{\frac{p\eta}{2}}s^{-\frac{p\eta}{2}}(1+|x|^p+|y|^p),
\end{eqnarray*}
where we use the fact that $1-e^{-x}\leq Cx^{\alpha\eta/2}$, for any $x>0$. The proof is complete.
\end{proof}

\begin{lemma}\label{L6.4}
For any $(x,y)\in H^{\gamma}\times H$ with $\gamma\in[0,1)$,  $\eta\in (0,1)$, $T>0$ and $p\in [1,\alpha)$, there exists a constant $C_{T}$ such that for any $\vare\in (0,1]$ and $t\in(0,T]$, we have
\begin{align}
\sup_{m\geq 1}\left[\mathbb{E}\|X^{m,\vare}_t\|^p_2\right]^{\frac{1}{p}}
\leq Ct^{-1+\frac{\gamma}{2}}\|x\|_{\gamma}+C_T\vare^{-\frac{\eta}{2}}(1+|x|+|y|). \label{Xvare2}
\end{align}
\end{lemma}
\begin{proof}
Note that for  any $t>0$, we have
\begin{eqnarray*}
X^{m,\vare}_t=\!\!\!\!\!\!\!\!&&e^{tA}x+\int_{0}^{t}e^{(t-s)A}B^m(X^{m,\vare}_{t},Y^{m,\vare}_{t})ds\nonumber\\
&&+\int_{0}^{t}e^{(t-s)A}\left[B^m(X^{m,\vare}_s,Y^{m,\vare}_{s})-B^m(X^{m,\vare}_{t},Y^{m,\vare}_{t})\right]ds+\int_{0}^{t}e^{(t-s)A}d\bar{L}^{m}_s  \nonumber\\
:=\!\!\!\!\!\!\!\!&&\sum^4_{i=1}I_i.
\end{eqnarray*}

For the term $I_{1}$, using \eqref{P3}, for any $\gamma\in (0,1)$ we have
\begin{align} \label{ABarX1}
\|e^{tA}x\|_2
\leq C t^{-1+\frac{\gamma}{2}}\|x\|_{\gamma}.
\end{align}

For the term $I_{2}$, we have
\begin{eqnarray*} \label{ABarX4}
\left[\mathbb{E}\|I_{2}\|^p_2\right]^{1/p}=\!\!\!\!\!\!\!\!&&\left(\mathbb{E}\left|(e^{tA}-I)B^m(X^{m,\vare}_{t},Y^{m,\vare}_{t})\right|^p\right)^{1/p}\\
\leq\!\!\!\!\!\!\!\!&&C\left[1+\left(\mathbb{E}|X^{m,\vare}_{t}|^p\right)^{1/p}+\left(\mathbb{E}|Y^{m,\vare}_{t}|^p\right)^{1/p}\right]\\
\leq\!\!\!\!\!\!\!\!&&C_T(1+|x|+|y| ).
\end{eqnarray*}

For the term $I_{3}$, using Minkowski's inequality and Lemma \ref{L6.3}, we obtain
\begin{eqnarray} \label{ABarX5}
\left[\mathbb{E} \| I_{3}\|^p_2\right]^{1/p}
\leq\!\!\!\!\!\!\!\!&&C\int_{0}^{t}\frac{1}{t-s}\left[\mathbb{E}\left|F^m_1(X^{m,\vare}_s,Y^{m,\vare}_{s})-F^m_1(X^{m,\vare}_{t},Y^{m,\vare}_{t})\right|^p\right]^{\frac{1}{p}}ds \nonumber\\
\leq\!\!\!\!\!\!\!\!&&
C\int_{0}^{t}\frac{1}{t-s}\left[\mathbb{E}\left|X^{m,\vare}_s-X^{m,\vare}_{t}\right|^p+\EE\left|Y^{m,\vare}_{s}-Y^{m,\vare}_{t}\right|^p\right]^{\frac{1}{p}}ds \nonumber\\
\leq\!\!\!\!\!\!\!\!&&C(1+|x|+|y|)\int_{0}^{t}\frac{1}{t-s}(t-s)^{\frac{\eta}{2}}s^{-\frac{\eta}{2}}ds\nonumber\\
&&+C(1+|x|+|y|)\int_{0}^{t}\frac{1}{t-s}\left(\frac{t-s}{\vare}\right)^{\frac{\eta}{2}}s^{-\frac{\eta}{2}}ds\nonumber\\
\leq\!\!\!\!\!\!\!\!&&C\vare^{-\eta/2}(1+|x|+|y|) .
\end{eqnarray}

For the term $I_{4}$, by \eref{LA} and assumption \ref{A2}, we easily have
\begin{eqnarray}
\left[\mathbb{E} \| I_{4}\|^p_2\right]^{1/p}\leq C_{p}\left(\sum_{k\in\mathbb{N}_{+}}\frac{\beta^{\alpha}_k}{\lambda^{1-\alpha}_k}\right)^{1/\alpha}\leq C_p.\label{ABarX6}
\end{eqnarray}

Combining \eref{ABarX1}-\eref{ABarX6} yields the desired result.
\end{proof}

\begin{remark}
By the same argument above, we can easily prove that for any $(x,y)\in H\times H$,  $T>0$ and $p\in [1,\alpha)$, there exists a constant $C_{T}$ such that for any $t\in(0,T]$, we have
\begin{align}
\sup_{m\geq 1}\left[\mathbb{E}\|\bar{X}^{m}_t\|^p_2\right]^{\frac{1}{p}}
\leq C_T t^{-1}(1+|x|). \label{barXvare2}
\end{align}
\end{remark}

\vspace{0.3cm}
Recall the  approximate equation \eref{Ga 1.3} to the averaged equation  \eref{1.3}.
Note that $\bar{X}^m$ is not the Galerkin approximation of $\bar{X}$, hence we have to check its approximation carefully.
\begin{lemma} \label{GA2}
For any $x\in H$, $T>0$ and $p\in (1,\alpha)$, we have
\begin{align}
\lim_{m\rightarrow \infty}\sup_{t\in[0,T]}\EE|\bar{X}^{m}_t-\bar{X}_t|^p=0. \label{FA2}
\end{align}
\end{lemma}
\begin{proof}
It is easy to see that for any $t>0$,
\begin{eqnarray*}
\bar{X}^{m}_t-\bar{X}_t=\!\!\!\!\!\!\!\!&&e^{tA}(x^m-x)+\int^{t}_{0}e^{(t-s)A}(\pi_m-I)\bar{B}(\bar{X}_s)ds\nonumber\\
&&+\int^{t}_{0}e^{(t-s)A}\left[\bar{B}^{m}(\bar{X}^{m}_s)-\pi_m \bar{B}(\bar{X}_s)\right]ds+\left[\int^t_0 e^{(t-s)A}d \bar L^{m}_s-\int^t_0 e^{(t-s)A}dL_s\right].
\end{eqnarray*}
Then for any $T>0$ and $p\in [1,\alpha)$, we have
\begin{eqnarray*}
&&\sup_{t\in [0,T]}\EE|\bar{X}^{m}_t-\bar{X}_t|^p\\
\leq\!\!\!\!\!\!\!\!&&C_p|x^m-x|^p+C_{p,T}\int^{T}_{0}\EE|(\pi_m-I)\bar{B}(\bar{X}_s)|^pds+C_{p,T}\int^{T}_{0}\EE|\bar{B}^{m}(\bar{X}^{m}_s)-\bar{B}^{m}(\bar{X}_s)|^pds\\
&&+C_{p,T}\EE\int^{T}_{0}|\bar{B}^{m}(\bar{X}_s)-\pi_m \bar{B}(\bar{X}_s)|^pds+C_p\sup_{t\in [0,T]}\EE\left|\int^t_0 e^{(t-s)A}d \bar L^{m}_s-\int^t_0 e^{(t-s)A}dL_s\right|^p\nonumber\\
\leq\!\!\!\!\!\!\!\!&&C_p|x^m-x|^p+C_{p,T}\int^{T}_{0}\EE|(\pi_m-I)\bar{B}(\bar{X}_s)|^pds+C_{p,T}\int^T_0 \EE|\bar{X}^{m}_t-\bar{X}_t|^pdt\\
&&+C_{p,T}\EE\int^{T}_{0}|\bar{B}^{m}(\bar{X}_s)-\pi_m \bar{B}(\bar{X}_s)|^pds+C_p\sup_{t\in [0,T]}\EE\left|\int^t_0 e^{(t-s)A}d \bar L^{m}_s-\int^t_0 e^{(t-s)A}dL_s\right|^p.
\end{eqnarray*}
By Gronwall's inequality, we get
\begin{eqnarray*}
&&\sup_{t\in [0,T]}\EE|\bar{X}^{m}_t-\bar{X}_t|^p\\
\leq\!\!\!\!\!\!\!\!&&C_{p,T}|x^m-x|^p+C_{p,T}\int^{T}_{0}\EE|(\pi_m-I)\bar{B}(\bar{X}_s)|^pds\nonumber\\
&&+C_{p,T}\EE\int^{T}_{0}|\bar{B}^{m}(\bar{X}_s)-\pi_m \bar{B}(\bar{X}_s)|^p ds+C_{p,T}\sup_{t\in [0,T]}\EE\left|\int^t_0 e^{(t-s)A}d \bar L^{m}_s-\int^t_0 e^{(t-s)A}dL_s\right|^p.
\end{eqnarray*}
By the a prior estimate of $\bar{X}$ and the dominated convergence theorem,
\begin{eqnarray}
&&\lim_{m\rightarrow \infty}\int^{T}_{0}\EE|(\pi_m-I)\bar{B}(\bar{X}_s)|^pds=0,\label{AT2.1}\\
&&\lim_{m\rightarrow \infty}\sup_{t\in [0,T]}\EE\left|\int^t_0 e^{(t-s)A}d \bar L^{m}_s-\int^t_0 e^{(t-s)A}dL_s\right|^p=0.\label{AT2.2}
\end{eqnarray}
Thus if we can prove for any $x\in H$,
\begin{eqnarray}
\lim_{m\rightarrow \infty}|\bar{B}^{m}(x)-\pi_m \bar{B}(x)|=0,\label{AT2.3}
\end{eqnarray}
Then by dominated convergence theorem, we get
\begin{eqnarray}
\lim_{m\rightarrow \infty}\EE\int^{T}_{0}|\bar{B}^{m}(\bar{X}_s)-\pi_m \bar{B}(\bar{X}_s)|^pds=0.\label{AT2.4}
\end{eqnarray}
Hence, \eref{FA2} holds by combining \eref{AT2.1}, \eref{AT2.2} and \eref{AT2.4}.

In fact, for any $t>0$,
\begin{eqnarray*}
\EE|Y^{x,0,m}_t-Y^{x,0}_t|\leq\!\!\!\!\!\!\!\!&&\int^{t}_{0}\EE|(\pi_m-I)F(x, Y^{x,0}_s)|ds+\int^{t}_{0}|F^{m}(x,Y^{x,0,m}_s)-F^m(x,Y^{x,0}_s)|ds \nonumber\\
&&+\EE\left|\int^t_0 e^{(t-s)A}d \bar Z^{m}_s-\int^t_0 e^{(t-s)A}dZ_s\right|\nonumber\\
\leq\!\!\!\!\!\!\!\!&&\int^{t}_{0}\EE|(\pi_m-I)F(x, Y^{x,0}_s)|ds+\EE\int^{t}_{0}|Y^{x,0,m}_s-Y^{x,0}_s|ds \nonumber\\
&&+\EE\left|\int^t_0 e^{(t-s)A}d \bar Z^{m}_s-\int^t_0 e^{(t-s)A}dZ_s\right|.
\end{eqnarray*}
By Gronwall's inequality,
\begin{eqnarray*}
\EE|Y^{x,0,m}_t-Y^{x,0}_t|
\leq\!\!\!\!\!\!\!\!&&e^{Ct}\left[\int^{t}_{0}\EE|(\pi_m-I)F(x, Y^{x,0}_s)|ds+\EE\left|\int^t_0 e^{(t-s)A}d \bar Z^{m}_s-\int^t_0 e^{(t-s)A}dZ_s\right|\right].
\end{eqnarray*}
As a consequence, it is easy to see
\begin{eqnarray*}
\lim_{m\rightarrow \infty}\EE|Y^{x,0,m}_t-Y^{x,0}_t|=0.
\end{eqnarray*}
By \eref{ergodicity1}, we have for any $t>0$,
\begin{eqnarray*}
|\bar{B}^{m}(x)-\pi_m \bar{B}(x)|\leq\!\!\!\!\!\!\!\!&&|\bar{B}^{m}(x)-\EE\pi_m B(x,Y^{x,0,m}_t)|+|\EE\pi_m B(x,Y^{x,0}_t)-\pi_m \bar{B}(x)| \nonumber\\
&&+|\EE\pi_m B(x,Y^{x,0,m}_t)-\EE\pi_m B(x,Y^{x,0}_t)|\nonumber\\
\leq\!\!\!\!\!\!\!\!&&Ce^{-\frac{(\lambda_1-L_F)t}{2}}+\EE|Y^{x,0,m}_t-Y^{x,0}_t|.
\end{eqnarray*}
Then by taking $m\rightarrow \infty$ firstly, then $t\rightarrow\infty$, we finally get \eref{AT2.3}. The proof is complete.

\end{proof}

\textbf{Acknowledgment}. This work is supported by the National Natural Science Foundation of China (11771187, 11931004, 12090011) and the Priority Academic Program Development of Jiangsu Higher Education Institutions.

\end{document}